\documentclass[11pt]{amsart}

\usepackage{times}
\usepackage[T1]{fontenc}
\usepackage{amssymb, amsthm, amsmath}
\usepackage{mathrsfs}

\usepackage[active]{srcltx}

\usepackage{todonotes}

\usepackage{setspace}
\usepackage{geometry}

\usepackage{hyperref}

\DeclareMathOperator{\acl}{acl}
\DeclareMathOperator{\dcl}{dcl}

\DeclareMathOperator{\stab}{Stab}

\DeclareMathOperator{\im}{im}

\DeclareMathOperator{\tp}{tp}\DeclareMathOperator{\stp}{stp}

\DeclareMathOperator{\mr}{RM}

\newtheorem{introtheorem}{Theorem}

\newtheorem{theorem}{Theorem}[section]

\newtheorem{claim}{Claim}[theorem]

\newtheorem{corollary}[theorem]{Corollary}

\newtheorem{fact}[theorem]{Fact}
\newtheorem{lemma}[theorem]{Lemma}

\newtheorem{proposition}[theorem]{Proposition}

\newtheorem*{gen-dif}{\fbox{{\large A}} \hypertarget{Agen-dif}{Gen-Dif}}

\newtheorem*{min-balln}{\fbox{{\large A}} \hypertarget{Amin-ball}{Cballs}}

\theoremstyle{definition}
\newtheorem{definition}[theorem]{Definition}
\newtheorem{example}[theorem]{Example}
\newtheorem{remark}[theorem]{Remark}
\newtheorem{question}[theorem]{Question}
\newtheorem{notation}[theorem]{Notation}
\newtheorem*{nota}{Notation}

\newcommand{\Cc}{{\mathbb{C}}}

\newcommand{\CN}{{\mathcal N}}
\newcommand{\CH}{{\mathcal H}}

\newcommand{\CM}{{\mathcal M}}

\newcommand{\CC}{{\mathcal C}}

\newcommand{\CG}{{\mathcal G}}

\newcommand{\CZ}{\mathcal Z}

\newcommand{\0}{\emptyset}

\renewcommand{\phi}{\varphi}

\newenvironment{claimproof}[1][\proofname]
  {%
    \proof[#1]%
  }
  {%
    \endproof%
  }

\newcommand{\dclindep}[1][]{%
  \mathrel{
    \mathop{
      \vcenter{
        \hbox{\oalign{\noalign{\kern-.3ex}\hfil$\vert$\hfil\cr
              \noalign{\kern-.7ex}
              $\smile$\cr\noalign{\kern-.3ex}}}
      }
    }\displaylimits_{#1}
  }
}

\long\def\symbolfootnote[#1]#2{\begingroup%
\def\thefootnote{\fnsymbol{footnote}}\footnote[#1]{#2}\endgroup}

\makeatletter

\def\Ind#1#2{#1\setbox0=\hbox{$#1x$}\kern\wd0\hbox to 0pt{\hss$#1\mid$\hss}
\lower.9\ht0\hbox to 0pt{\hss$#1\smile$\hss}\kern\wd0}
\def\ind{\mathop{\mathpalette\Ind{}}}
\def\Notind#1#2{#1\setbox0=\hbox{$#1x$}\kern\wd0\hbox to 0pt{\mathchardef
\nn=12854\hss$#1\nn$\kern1.4\wd0\hss}\hbox to
0pt{\hss$#1\mid$\hss}\lower.9\ht0 \hbox to
0pt{\hss$#1\smile$\hss}\kern\wd0}

\def\sub{\subseteq}

\makeatother

\title{Very Ampleness in Strongly Minimal Sets}

\date{\today}
\author{Benjamin Castle}
\address{Department of Mathematics, Ben Gurion University of the Negev, Be'er-Sheva 84105, Israel}
\email{bcastle@math.berkeley.edu }

\author{Assaf Hasson}
\address{Department of Mathematics, Ben Gurion University of the Negev, Be'er-Sheva 84105, Israel}
\email{hassonas@math.bgu.ac.il}

 \thanks{Both authors were partially supported by ISF grant No. 555/21. The second Author was supported by Ben Gurion Kreitman fellowship.}

\begin{document}

\begin{abstract}

   Inspired by very ampleness of Zariski Geometries, we introduce and study the notion of a very ample family of plane curves in any strongly minimal set, and the corresponding notion of a very ample strongly minimal set (characterized by the definability of such a family). We show various basic properties; for example, any strongly minimal set internal to an expansion of an algebraically closed field is very ample, and any very ample strongly minimal set non-orthogonal to a strongly minimal set $Y$ is internal to $Y$. We then apply these results with Zilber's restricted trichotomy to characterize using very ampleness those structures $\mathcal M=(M,\dots)$ interpreted in an algebraically closed field which recover all constructible subsets of powers of $M$. Next we show that very ample strongly minimal sets admit very ample families of plane curves of all dimensions, and use this to characterize very ampleness in terms of definable pseudoplanes. Finally, we show that divisible strongly minimal groups are very ample, and deduce -- answering an old question of G. Martin, \cite{Martin} -- that in a pure algebraically closed field, $K$ there are no reducts between $(K,+,\cdot)$ and $(K, \cdot)$. 
\end{abstract}

\maketitle


\tableofcontents

\section{Introduction}

Many ideas in stability theory (such as Zilber's classification of totally categorical structures, see e.g. \cite{ZilBook}) build heavily on or are inspired by Zilber's weak trichotomy theorem, stating that any strongly minimal set is either locally modular or defines a rank 2 pseudoplane (recall that a \textit{pseudoplane} is an abstract incidence relation $I\subset P\times L$ between `points' and `lines' such that each point lies on infinitely many lines, each line contains infinitely many points, any two lines have finite intersection, and any two points lie on only finitely many common lines). Indeed, Zilber's (strong) trichotomy conjecture, in one of its many variants, asserted that every uncountably categorical pseudoplane should be mutually interpretable with an algebraically closed field. 

While Zilber's conjecture is false (\cite{Hr1}), it has been shown to hold in many instances. For example, the first author has recently shown that the conjecture holds assuming the pseudoplane is itself interpreted in an algebraically closed field of characteristic zero (\cite{CasACF0}), and the second author and D. Sustretov have shown the same for rank preserving interpretations in all characteristics (\cite{HaSu}) (this uses the well known fact that  uncountably categorical pseudoplanes are not 1-based).

Note, however, that though a non-locally modular strongly minimal set, $S$, always interprets a pseudoplane, the set of `points' of the pseudoplane might not be identifiable with $S^2$. In general, the points could be imaginary, which may complicate the process of recovering a field. Thus, for example, it seems that Rabinovich's work \cite{Ra} would be considerably simpler if one assumes the set of points of the given pseudoplane to be 
the affine plane\footnote{See discussion on p.2-3 of \cite{Ra}. Unfortunately, the manuscript [9] referred to in that text does not seem to be available.}. In practice, in the proof of \cite{CasACF0}, as well as in several other instances where the trichotomy was established, it turned out to be more useful to work directly with a rank 2 family of plane curves in $S$, even if such a family does not form a (combinatorially cleaner) pseudoplane.

In fact, we do not even know whether the interpreted pseudoplane (provided by Zilber in a non-locally modular strongly minimal set, $D$) can be constructed on the square of \textit{any} strongly minimal set (let alone $D$). However, we do know by Zilber's thesis (Lemma III.2.2, according to Hodges' introduction to \cite{Ra}) that an uncountably categorical structure interpreting some pseudoplane must interpret a two-dimensional faithful family of plane curves in every strongly minimal set. The crucial difference here is that such a family of plane curves allows certain exceptional pairs of points (called \textit{semi-indistinguishable pairs} in \cite{CasACF0}) to lie on infinitely many common curves.

The effect has been that working directly with families of plane curves allows for more direct access to a background geometry, in particular recovering (in reasonably geometric settings) approximations within a strongly minimal structure of such notions as ``closure," ``frontier," and ultimately ``tangency" (leading to the recovery of a field); however the combinatorial aspects of the construction become quite technical, essentially due to the presence of semi-indistinguishable points. These points tend to prohibit the `full' recovery of the underlying geometry -- indeed, already in \cite{Martin} it was pointed out (in an example attributed to Hrushovski -- see Example \ref{canonical example ex}) that a reduct of an algebraically closed field need not be constructible in any copy of the field it interprets (or even $\CM$-definably isomorphic to such a constructible set).

In their seminal work, \cite{HrZil}, Hrushovski and Zilber identify a condition they call \emph{very ampleness}, assuring that a Zariski Geometry, $\CZ$,  is algebraic in the above sense, i.e., that not only does it interpret a field, $K$, but it is definably isomorphic (as a Zariski Geometry) to a smooth $K$-algebraic curve. Very ampleness (as defined by Hrushovksi and Zilber) does not go so far as to require a pseudoplane on the set $Z^2$, but it is a related notion -- namely, it asserts the existence of a family of plane curves that separates points in $Z^2$ in an obvious sense.

In the present paper, we introduce a similar notion of very ampleness for arbitrary strongly minimal sets, characterized (Corollary \ref{very ample iff very ample family cor}) by the presence of families of plane curves with no semi-indistinguishable pairs. In fact we show that the presence of one such family implies the existence of such families of arbitrary dimension (Corollary \ref{C: families of all dims}); using this, we conclude (Proposition \ref{P: VA = nice PP}) that a strongly minimal set $D$ is very ample if and only if there is a definable pseudoplane whose points are (a generic subset of) $D^2$. Moreover, we show (Proposition \ref{P: zg va equiv}) that in the context of Zariski Geometries, our version of very ampleness is a necessary and sufficient condition for a Zariski geometry to be algebraic, and is thus equivalent to the original notion. 

We also apply very ampleness to the study of structures interpretable in algebraically closed fields. We show (Corollary \ref{C: field va}) that every algebraic curve over such a field is very ample. We then show (Proposition \ref{internality prop}) that a very ample strongly minimal set is internal to every strongly minimal set it interprets -- a phenomenon allowing for the identification of very ample ACF-interpreted structures with algebraic curves. In particular (noting the known instances of the Zilber trichotomy in the relevant structures -- see \cite{CasACF0} and \cite{HaSu}) we conclude the following analogue of the main theorem on Zariski geometries (see Theorem \ref{very ample acf char thm}):

\begin{introtheorem}\label{T: intro 0} Let $K$ be an algebraically closed field, let $M$ be constructible over $K$, and let $\mathcal M=(M,\dots)$ be a strongly minimal reduct of the full $K$-induced structure on $M$. Assume $\mathcal M$ satisfies the Zilber trichotomy. Then the following are equivalent:
\begin{enumerate}
    \item $\mathcal M$ is very ample.
    \item $\mathcal M$ is isomorphic, outside a finite set, to the full $K$-induced structure on some irreducible algebraic curve over $K$.
    \item Every constructible subset of every power of $M$ is definable in $\mathcal M$.
\end{enumerate}
\end{introtheorem}

For arbitrary ACF-interpreted structures, we get the following similar result (Theorem \ref{T: higher relics}):

\begin{introtheorem}\label{T: intro 1} Let $K$ be an algebraically closed field, and assume the Zilber trichotomy holds for strongly minimal structures interpreted in $K$. Let $M$ be constructible over $K$, and let $\mathcal M=(M,\dots)$ be an arbitrary reduct of the $K$-induced structure on $M$. Then the following are equivalent:
\begin{enumerate}
    \item $\mathcal M$ is almost strongly minimal, and every strongly minimal set in $\mathcal M$ is very ample.
    \item Every constructible subset of every power of $M$ is definable in $\mathcal M$.
\end{enumerate}
\end{introtheorem}

We then apply Theorem \ref{T: intro 0} in the case of groups, recovering in particular the following (Theorem \ref{group trichotomy thm}): 

\begin{introtheorem}
    Let $K$ be an algebraically closed field, and let $(G,\cdot)$ be a $1$-dimensional divisible algebraic group over $K$. Let $\mathcal G^{\mathrm{Zar}}$ be the full $K$-induced structure on $G$, and let $\mathcal G^{\mathrm{lin}}$ be the structure endowing $G$ with the group operation and all $K$-definable endomorphisms of $G$. Then there are no intermediate structures between $\mathcal G^{\mathrm{lin}}$ and $\mathcal G^{\mathrm{Zar}}$.
\end{introtheorem}

This expands the main result of \cite{MaPi} and can be seen as an algebraically closed field analogue of a recent similar result of Abu Saleh and Peterzil \cite{aSaPet} for real closed fields. 

Along the way we prove several technical results that may be of interest on their own right. We single out the following (Theorem \ref{T: nlm implies non affine}):
\begin{introtheorem}\label{T: intro 2}
    If $\CG$ is a strongly minimal expansion of a group, then $\CG$ is not locally modular if and only if there exists a definable $X\sub G^2$ that is not a finite boolean combination of cosets of definable subgroups of $G$. 
\end{introtheorem}

This result, while well known among experts, does not seem to exist in writing in full generality, and since we needed it in the present work we took the opportunity to give the details. \\

The paper is written with an eye toward the proofs of Theorems \ref{T: intro 0}-\ref{T: intro 2} above, but also as a possible reference for future work around Zilber's Trichotomy and other questions related to the fine structure of strongly minimal sets and the structures they interpret. For that reason, some of the proofs are stated and proved in somewhat greater generality than is actually needed. Thus, for example, the existence of arbitrarily large  very ample families of plane curves (Proposition \ref{P: going up}) is not explicitly used in the text, but significant parts of its proof (e.g. very ampleness of algebraic curves, Corollary \ref{C: field va}) are essential for our arguments. \\

\noindent \emph{Acknowledgements} We would like to thank E. Hrhoshovski for clarifications regarding the definition of very ampleness in \cite{HrZil}. 

\section{Notation and preliminaries}
Throughout, we work in a saturated enough model of a stable theory. Except for the general definition of very ample types (invoked only in Corollary \ref{C: 1-based iff no va type}) the work could be carried out in any structure where all definable strongly minimal sets are stably embedded. Since stable embeddedness is equivalent to uniform stable embeddedness, given a strongly minimal set $X$ there is no harm in assuming that all definable families of subsets of $X^n$ we are considering are parameterised by $X$-definable sets. 

We are using standard model theoretic terminology. For basic concepts such as Morley rank, strong minimality, local modularity, 1-based theories etc. we refer the reader to any textbook covering the first chapters of geometric stability theory such as \cite{MaBook}, \cite{PillayBook} or \cite{TZ}. 

After fixing a structure $\mathcal M$, unless otherwise stated, the word \textit{definable} refers to definability in $\mathcal M^{\textrm{eq}}$ with parameters. All parameter sets are smaller than the level of saturation of $\mathcal M$, and are usually denoted, $A$, $B$. We use the standard model theoretic abuse of notation and write $a\in M$ instead of $a\in M^{|a|}$, allowing also $a$ to be an element of an imaginary sort. 

If $X$ is a strongly minimal definable set, the notation $X^{\textrm{eq}}$ refers to the union of all sorts formed by quotienting definable subsets of powers of $X$ by definable equivalence relations. A definable set in $X^{\textrm{eq}}$ will be called stationary if it has a unique generic type over any set of parameters defining it. It follows by uniform stable embeddedness that every stationary set $S$ in $X^{\textrm{eq}}$ has a canonical base in $X^{\textrm{eq}}$ (i.e. the canonical base of its generic type), which we denote $\operatorname{Cb}(S)$. We use $\dim$ and \textit{dimension} to denote Morley rank of definable sets in $X^{\textrm{eq}}$. If $Y\sub Z$ are definable sets in $X^{\textrm{eq}}$, we say that $S$ is large in $Z$ if $\dim(Z\setminus Y)<\dim(Z)$ and $Y$ is small in $Z$ if $\dim(Y)<\dim(Z)$. A generic subset of $Z$ is a definable subset of $Z$ that is not small. Thus, if $Y\subset Z$ and $Z$ is stationary, $Y$ is generic in $Z$ if and only if it is large in $Z$.

We will use implicitly and without further reference the definability of Morley rank in $X^{\textrm{eq}}$ for strongly minimal $X$ (see \cite{BalUCCT}). We will also use the resulting fact that, given a small parameter set $A$, a definable set $C$ in $X^{\textrm{eq}}$ is always a generic member of an $A$-definable family of $\CM$-definable subsets of the same Morley rank. 
\\

By a \textit{curve} (in $X^{\textrm{eq}}$) we mean a $1$-dimensional definable set, and by a \textit{plane curve} in $X$, we mean a curve in $X^2$. Note that we do not require curves to be strongly minimal. A plane curve $C$ is \textit{trivial} if one of its projections has an infinite fibre. If $C, D$ are non-trivial plane curves, then their composition $D\circ C$ is -- in analogy with the composition of functions -- the curve $\{(x,z)): (\exists y)(x,y)\in C \land (y,z)\in D\}$. 

A definable family of plane curves $\CC:=\{C_t: t\in T\}$ is \emph{faithful} if $t\neq t'$ implies $\dim(C_t\vartriangle C_t')=0$. More generally, a definable family of definable sets in $X^{\textrm{eq}}$ is faithful if the symmetric difference of any two members of the family is small in both. The total space (or the graph) of the family $\CC$ is denoted $C:=\{(x,t)\in 
X^2\times T: x\in C_t\}$.

A definable set $C$ is \emph{almost contained} in a definable set $C'$ if $\dim(C'\setminus C))<\dim(C)$. The set $C$ is \emph{almost equal} to $C'$ if $C$ is almost contained in $C'$ and $C'$ is almost contained in $C$ (equivalently, if $\dim(C\vartriangle C')<\dim(C)=\dim(C')$). This is a definable equivalence relation (on families of definable sets). It is now a standard and easy exercise using imaginary elements to show that, up to a small correction, any stationary definable set $C$ in $X^{\textrm{eq}}$ (for $X$ strongly minimal and $\emptyset$-definable) is a generic member of a $\0$-definable \emph{faithful} family of $\dim(C)$-dimensional sets. 

A strongly minimal set $X$ is locally modular, if any definable faithful family of plane curves in $X^2$ is (at most) 1-dimensional. By this we mean that if such a family is parameterised by a set $T$ then $\dim(T)$,  the dimension of the family, is at most $1$.  We will systematically use (without further reference) the equivalence of local modularity and 1-basedness, which -- in turn -- is equivalent to the fact that $\dim(\mathrm{Cb(S)})\le 1$ for any strongly minimal $S\sub X^2$ (\cite[\S II Proposition 2.6]{PillayBook}).

\section{Very Ampleness: First Properties}

\subsection{Definition and first examples}

We start with a general definition, that we later on specialise almost exclusively to the strongly minimal setting:

\begin{definition}\label{very ample type def} Let $p$ be a stationary type, with canonical base $c$, and let $A$ be a set of parameters. We say that $p$ is \textit{very ample over} $A$ if for any two distinct realizations $x,y\models p$, $x$ and $c$ fork over $Ay$.
\end{definition}
\begin{definition}\label{very ample sm plane curve def}
Let $X$ be a strongly minimal set, definable over a set $A$. Let $C\subset X^2$ be a strongly minimal plane curve in $X$. We say that $C$ is \textit{very ample in} $X$ \textit{over} $A$ if the generic type of $C$ is very ample over $A$.
\end{definition}
\begin{definition} Let $X$ be a strongly minimal set. We say that $X$ is \textit{very ample} if, for some set $A$ such that $X$ is $A$-definable, there is a very ample strongly minimal plane curve in $X$ over $A$.
\end{definition}

\begin{definition}
Let $\mathcal M$ be a strongly minimal structure. We say that $\mathcal M$ is very ample if its universe is very ample as a strongly minimal set.
\end{definition}

\begin{example}\label{acf line ex} Let $(K,+,\cdot)$ be an algebraically closed field, and let $(a,b)\in K^2$ be generic. We claim that the line $L_{(a,b)}$ defined by $y=ax+b$ is very ample in $K$ over $\emptyset$, and thus $K$ is very ample. To see this, first note that the canonical base of the generic type of $L_{(a,b)}$ is just $(a,b)$. Now let $z_1=(x_1,y_2)$ and $z_2=(x_2,y_2)\in L_{(a,b)}$ be distinct generics over $(a,b)$. One then easily computes that $\dim(ab/z_2)=1$, essentially because the family of lines through a point in the plane is one-dimensional. On the other hand, because any two distinct points determine exactly one line, we have $\dim(ab/z_1z_2)=0$, which gives the desired dependence.
\end{example}

\begin{example}\label{canonical example ex} Suppose $(K,+,\cdot)$ is an algebraically closed field, and consider the reduct $\mathcal M$ of $K$ obtained by endowing $K$ with the relations $x^2+y^2=z^2$ and $x^2y^2=z^2$. Then $\mathcal M$ is non-locally modular, but is not very ample. Indeed, suppose $C\subset K^2$ is a strongly minimal plane curve which is very ample over $A$, and let $c=\operatorname{Cb}(C)$. Let $(x,y)\in C$ be generic over $Ac$. Without loss of generality, assume $x$ is generic in $K$ over $Ac$ (otherwise $y$ is). Note that any function $\sigma: K\to K$ with the property that $\sigma(x)\in \{x,-x\}$ for all $x$  is an automorphism of $\mathcal M$. In particular, there is an automorphism fixing $Ac$ pointwise but sending $x$ to $-x$, and $y$ to some $z=\pm y$. It follows that $\tp(-x,z/Ac)=\tp(x,y/Ac)$, so $(-x,z)$ is generic in $C$ over $Ac$. Then, by very ampleness, $(-x,z)$ should fork with $c$ over $Axy$; but this is impossible because $(-x,z)$ is algebraic over $(x,y)$. 
\end{example}

\begin{remark}
    We note that there is a significant difference between a strongly minimal set, $X$, being very ample (meaning that it admits a very ample plane curve) and its generic type being very ample. As we will see in the sequel, the existence of a (strongly minimal) very ample complete type is equivalent to the structure being non 1-based. The former condition is stronger and is equivalent to the existence of a definable pseudoplane whose set of points is generic in $X^2$. Except for Lemma \ref{C: 1-based iff no va type} the latter property is reserved for plane curves in $X$ (and we will say that $C$ is very ample in $X$), so hopefully no confusion can arise. 
\end{remark}

\subsection{Very Ampleness and Non-local Modularity}
It follows directly from the definition that a very ample strongly minimal set is not locally modular. We give the details:

\begin{lemma}\label{L: very ample implies nlm} Let $X$ be strongly minimal and definable over $A$, let $C$ be a strongly minimal plane curve in $X$, and let $c=\operatorname{Cb}(C)$. If $C$ is very ample in $X$ over $A$, then $\dim(c/A)\geq 2$. In particular, if $X$ is very ample then $X$ is not locally modular.
\end{lemma}
\begin{proof}
Let $x$ and $y$ be independent generics in $C$ over $Ac$. By very ampleness, $c$ forks with $x$ over $Ay$, which implies in particular that $\dim(c/Ay)\geq 1$.
\begin{claim} $\dim(y/A)=2$.
\end{claim}
\begin{proof} Otherwise, $y$ would belong to a one-dimensional set over $A$, which must have $C$ as a strongly minimal component. But this would force $c\in\acl(A)$, contradicting that $\dim(c/Ay)\geq 1$.
\end{proof}
Now by the claim, and since $\dim(y/Ac)=1$ by definition, it follows that $y$ forks with $c$ over $A$. By symmetry this implies $$\dim(c/A)>\dim(c/Ay)\geq 1,$$ i.e. $\dim(c/A)\geq 2$.
\end{proof}

We will see later on, Corollary \ref{C: 1-based iff no va type}, that the weaker condition of $X^{\textrm{eq}}$ containing a very ample type is equivalent to non-local modularity. 

\begin{remark}\label{R: non trivial} It also follows from this lemma that if a strongly minimal plane curve $C$ is very ample in $X$ over some set, then $C$ is \textit{non-trivial}, i.e. both projections $C\rightarrow X$ are finite-to-one. Indeed, otherwise $C$ would agree up to finitely many points with $\{c\}\times X$ or $X\times\{c\}$ for some $c\in X$, implying that $c=\operatorname{Cb}(X)$ has dimension at most 1 over $A$.
\end{remark}

Note that if $X,Y$ are strongly minimal sets such that $Y\vartriangle X$ is finite, then for any non-trivial strongly minimal plane curve $C\sub X^2$ the set $C\setminus Y^2$ is finite, so $C\cap Y^2$ and $C$ have the same generic type. In particular, if $C$ is very ample in $X$ then $C\cap Y^2$ is very ample in $Y$. Thus $X$ is very ample if and only if $Y$ is.

Before moving on we mention the following useful and related observation, which may help to further clarify the definition. Roughly, it shows that the non-very ampleness of a non-locally modular strongly minimal set always arises from an interalgebraic relation on the plane:

\begin{lemma}\label{L: CE are inter-alg}
    Let $X$ be strongly minimal and $A$-definable. Let $C$ be a plane curve in $X$, definable over $B\supset A$. Assume that for every strongly minimal component $S$ of $C$ we have $\dim(\operatorname{Cb}(S)/A)\geq 2$. Let $x,y\in C$ be any two distinct generics over a parameter set $B\supseteq A$. Then either $x$ forks with $B$ over $Ay$, or $x$ and $y$ are interalgebraic over $A$.
\end{lemma}
\begin{proof}
   Throughout we assume $A=\emptyset$. Let $S\subset C$ be a strongly minimal component containing $x$ and let $s=\mathrm{Cb}(S)$. Note that $s\in\acl(B)$ and $\dim(x/s)=1$ (since $x$ is generic in $C$). Now, by assumption $\dim(s)\geq 2$; in particular $s\notin\acl(\emptyset)$. So $\tp(x/s)$ must fork over $\emptyset$, and thus $\dim(x)=2$. We now argue in two cases:
   \begin{itemize}
       \item First suppose $x\in\acl(y)$. Since $\dim(x)=2$, this implies $\dim(y)=2$, and thus that $x$ and $y$ are interalgebraic, proving the lemma in this case.
       \item Now suppose $x\notin\acl(y)$; we show that $x$ forks with $B$ over $y$. By assumption we have $\dim(x/y)\geq 1$. Note also that $\dim(x/By)\leq 1$ (as $x\in C_t$). Now suppose toward a contradiction that $x\ind_yB$. Then clearly we have $\dim(x/y)=\dim(x/By)=1$. Since $s\in\acl(B)$, also $\dim(x/Bsy)=1$. In particular $\tp(x/Bsy)$ (i.e. the generic type of $S$ over $Bsy$) does not fork over $y$; so since $s=\operatorname{Cb}(S)$, we get that $s\in\acl(y)$. Since $y\in X^2$, this gives $\dim(s)\leq 2$. Combined with the assumption that $\dim(s)\leq 2$, we therefore have $\dim(s)=2$, which implies that $s$ and $y$ are interalgebraic. But then $\dim(y/s)=\dim(y/Bs)=0$, and since $s\in\acl(B)$ this implies $\dim(y/B)=0$, contradicting that $y$ is generic in $C$ over $B$.
    \end{itemize}
\end{proof}

\subsection{Extending the Definition}
In many applications, it is more convenient to work with families of plane curves that are not necessarily generically strongly minimal. We extend our definition of very ampleness to such (families of) curves.  To do this, we observe that the choice of defining parameter for a plane curve, $C$, does not affect Definition \ref{very ample sm plane curve def} in any significant way:

\begin{lemma}\label{L: Ind on parms} Let $X$ be strongly minimal and $A$-definable, and let $C$ be an $At$-definable plane curve in $X$ for some tuple $t$. The following are equivalent:
\begin{enumerate}
    \item For any two distinct $x, y\in C$ each generic over $At$, $x$ forks with $t$ over $Ay$.
    \item For any two distinct $x, y\in C$ each generic over $At$, $\dim(t/Axy)\leq\dim(t/A)-2$.
    \item For every strongly minimal set $S\subset C$, $\dim(\operatorname{Cb}(S)/A)\geq 2$, and for any two distinct $x, y\in C$ each generic over $At$, $x$ and $y$ are not interalgebraic over $A$.
\end{enumerate}
\end{lemma}
\begin{proof} Throughout, we assume $A=\emptyset$. 

(1) $\Rightarrow$ (2): Assume (1), and let $x$ and $y$ be distinct generics of $C_t$ over $t$.
\begin{claim} $\dim(y)=2$.
\end{claim}
\begin{proof}
Let $y'$ be an independent realization of $\tp(y/t)$ over $ty$. Then by (1), $y$ forks with $t$ over $y'$. Since $\dim(y/t)=1$ by assumption, we obtain $\dim(y/y')=2$, which implies the claim.
\end{proof}
Now since $\dim(y)=2$ and $\dim(y/t)=1$, it follows that $\dim(t/y)=\dim(t)-1$. By (1) again, $t$ forks with $x$ over $y$, which gives $\dim(t/xy)\leq\dim(t)-2$, as desired.

(2) $\Rightarrow$ (3): Assume (2). First let $S$ be a strongly minimal component of $C$, with canonical base $s$. Let $x$ and $y$ be independent generics in $S$ over $st$. Then $\dim(txy)=\dim(t)+2$, and by (2) $\dim(t/xy)\leq\dim(t)-2$, which combined gives $\dim(xy)=4$. On the other hand $\dim(sxy)=\dim(s)+2$, so since $\dim(xy)=4$ we get $\dim(s/xy)=\dim(s)-2$. In particular $\dim(s)-2\geq 0$, so $\dim(s)\geq 2$.

Now to complete the proof of (3), suppose $x$ and $y$ are any two distinct generics of $C$ over $t$. Then $\dim(txy)\geq\dim(tx)=\dim(t)+1$, while by (2) $\dim(t/xy)\leq\dim(t)-2$. It follows that $\dim(xy)\geq 3$, while if $x$ and $y$ were interalgebraic we would have $\dim(xy)\leq 2$. Thus, we have shown (3).

(3) $\Rightarrow$ (1): Assume (3). Let $x$ and $y$ be a counterexample to (1). By (3), the hypotheses of Lemma \ref{L: CE are inter-alg} are satisfied, and applying the lemma gives that $x$ and $y$ are interalgebraic over $A$, contradicting (3). 
\end{proof}

\begin{corollary}
Let $X$ be strongly minimal and $A$-definable. For $i=1,2$ let $C_i$ be an $At_i$-definable plane curve in $X$ for some tuple $t_i$. Assume that $C_1$ and $C_2$ have finite symmetric difference. Then the equivalent conditions of Lemma \ref{L: Ind on parms} hold for $C_1$ over $At_1$ if and only if they hold for $C_2$ over $At_2$.
\end{corollary}
\begin{proof}
We use (3) of the lemma. Since the $C_i$ have finite symmetric difference, the canonical bases of their strongly minimal components coincide. So, by symmetry, it suffices to show that the failure of the second clause of (3) for $C_1$ and $t_1$ implies the failure of the same clause for $C_2$ and $t_2$. To that end, let $x$ and $y$ be distinct generics in $C_1$ over $At_1$, and assume they are interalgebraic over $A$. Let $(x',y')$ be an independent realization of $\tp(xy/At_1)$ over $At_1t_2$. Then $x'$ and $y'$ are independent generics over $C_2$ over $At_2$, and are also interalgebraic over $A$, as desired.
\end{proof}

The following is now well-defined:

\begin{definition}\label{very ample general plane curve def} Let $X$ be strongly minimal and $A$-definable. Let $C$ be a plane curve in $X$. Then $C$ is \textit{very ample in} $X$ \textit{over} $A$ if any of the equivalent conditions of Lemma \ref{L: Ind on parms} hold for $C$ over $At$, for some (equivalently any) tuple $t$ such that $C$ is $At$-definable.
\end{definition}

It is immediate from the last corollary that Definitions \ref{very ample sm plane curve def} and \ref{very ample general plane curve def} coincide for strongly minimal plane curves, since up to a finite set any such curve is definable over its canonical base. Let us also point out the following: 

\begin{corollary}\label{C: VA components} Let $X$ be strongly minimal and definable over $A$. Let $C$ be a plane curve in $X$. If $C$ is very ample in $X$ over $A$, then so is every strongly minimal component of $C$.
\end{corollary}

\begin{proof}
Assume $C$, and all of its strongly minimal components, are definable over $At$ for some tuple $t$, and let $S\subset C$ be any such component. Let $x\neq y$ be generics in $S$ over $At$. Then $x$ and $y$ are also generics of $C$ over $At$; so by the very ampleness of $C$ over $A$, $x$ forks with $t$ over $Ay$, which shows that $S$ is also very ample over $A$.
\end{proof}

The converse of the last corollary is not true: there can be non very ample plane curves every component of which is very ample, as in the following example:
\begin{example}  Let $\mathcal G=(G,+,\dots)$ be a strongly minimal expansion of a group, and let $S\subset G^2$ be a $\emptyset$-definable strongly minimal plane curve whose generic type has trivial stabilizer (for example $\mathcal G$ could be the additive group of the complex field, and $S$ could be the graph of $y=x^2$). Now let $a\in G^2$ be generic, let $S_a:=S+a$, and let $C=S_a\cup -S_a$. It is easy to see that for $s\in S_a$ generic over $a$, the pair $s,-s\in C$ contradicts the very ampleness of $C$ (by (3) of Lemma \ref{L: Ind on parms}). However, each of $S_a$ and $-S_a$ is very ample (see Example \ref{group ex} below). 
\end{example}

\subsection{Very Ampleness in Families}

The reader may find our definition of very ampleness foreign in light of the analogous notion  in \cite{HrZil} (a precise statement on the relation of our notion of very ampleness and the original term in the context of Zariski Geometries requires additional preparation, and can be found in Section \ref{Ss: ZG}). We hope that the present subsection will help to explain for now why our definition captures the original idea of the term. To that end, we now define a notion of very ampleness for families of plane curves:

\begin{nota}
If $S$ is definable, and $\{X_t:t\in T\}$ is a definable family of subsets of $S$, then for $s\in S$ we denote the set $\{t\in T:s\in X_t\}$ by $X^s$.
\end{nota}

\begin{definition}
Let $X$ be strongly minimal and $A$-definable, and let $\mathcal C=\{C_t:t\in T\}$ be an $A$-definable family of plane curves in $X$.
\begin{enumerate}
    \item $\mathcal C$ is \textit{very ample} if for any $x\neq y\in M^2$ we have $\dim(C^x\cap C^y)\leq\dim T-2$.
    \item $\mathcal C$ is \textit{generically very ample} if for any generic $t\in T$ over $A$ and any distinct $x,y\in C_t$ generic over $At$ we have $\dim(t/Axy)\leq\dim T-2$. 
\end{enumerate}
\end{definition}

Observe that, in particular, a generically very ample family of plane curves is at least two-dimensional. 

\begin{example}\label{acf lines ex} Let $(K,+,\cdot)$ be an algebraically closed field. Then, similarly to Example \ref{acf line ex}, one shows easily that the family of lines $y=ax+b$ is very ample, because any two distinct points determine exactly one line.
\end{example}

The following example is well known and easy: 

\begin{example}\label{group ex} Let $(G,+)$ be a strongly minimal expansion of a group, and let $C\subset G^2$ be a strongly minimal plane curve. Then the family of translates of $C$, $\{C+t:t\in G^2\}$, is very ample if and only if the generic type of $C$ has trivial stabilizer.
\end{example}

Proposition \ref{P: family VA iff curve VA} below, in addition to the ensuing two useful corollaries, shows the relationship between the notions of very ampleness we have defined so far (of a plane curve, of a strongly minimal set, and of a family of plane curves).

\begin{proposition}\label{P: family VA iff curve VA}
Let $X$ be strongly minimal and $A$-definable, let $\mathcal C=\{C_t:t\in T\}$ be an $A$-definable family of plane curves in $X$, and let $C\subset X^2\times T$ be the graph of $\mathcal C$. Then the following are equivalent:
\begin{enumerate}
    \item $\mathcal C$ is generically very ample.
    \item There is a non-generic $A$-definable set $Z\subset C$ such that $C-Z$ is the graph of a very ample family of plane curves in $X$.
    \item For any generic $t\in T$ over $A$, $C_t$ is very ample in $X$ over $A$.
\end{enumerate}
\end{proposition}
\begin{proof}
The equivalence of (1) and (3) is Lemma \ref{L: Ind on parms}. 

(1) $\Rightarrow$ (2): Let $p(x,y,z)$ be a partial type over $A$ asserting that $x\neq y\in X^2$ and that $(x,z)$ and $(y,z)$ are both generics in $C$ over $A$. It is a restatement of (1) that for $(x,y,z)\models p$ we have $\dim(z/Axy)\leq\dim(T)-2$. By compactness, this is witnessed by a finite part of $p$. We thus obtain a non-generic $Z\subset C$ such that for $x\neq y$ we have $\dim((C-Z)^x\cap(C-Z)^y)\leq\dim T-2$. This almost implies (2); technically, though, one should add to $Z$ all $C_t$ having only finitely many points surviving in $C-Z$ (so that $C-Z$ is a family of plane curves); one should then verify that $Z$ is still non-generic in $C$, and that set of `surviving' $t\in T$ in $C-Z$ is generic in $T$. Both of these are easy to do. 

(2) $\Rightarrow$ (1): Let $Z$ be as in (2), and let $p(x,y,z)$ be the same type considered in the case above. We show that for $(x,y,z)\models p$ we have $\dim(z/Axy)\leq\dim T-2$, which as above is equivalent to (1). So take such $(x,y,z)$. We may assume that $(x,y,z)\ind_AB$ where $Z$ is $B$-definable; thus $(x,y,z)$ is generic in $C$ over $AB$, so $(x,y,z)\in C-Z$. Then by very ampleness, $\dim(z/ABxy)\leq\dim(T)-2$, and again since $(x,y,z)\ind_AB$ this implies $\dim(z/Axy)\leq 2$, as desired.

\end{proof}

To sum up, we show that a plane curve is very ample precisely when it coincides,  up to a finite set, with a generic member of a very ample family: 

\begin{corollary}\label{C: curves to families} Let $X$ be strongly minimal and $A$-definable. Let $C$ be a plane curve in $X$. Then the following are equivalent:
\begin{enumerate}
    \item $C$ is very ample in $X$ over $A$.
    \item $C$ has finite symmetric difference with an $A$-generic member of an $A$-definable very ample family of plane curves in $X$.
    \item Whenever $T$ is a stationary definable set, $\mathcal C=\{C_t:t\in T\}$ is an $\acl(A)$-definable family of plane curves in $X$, and $t\in T$ is a generic element over $\acl(A)$ such that $C$ and $C_t$ have finite symmetric difference, the family $\mathcal C$ is generically very ample.
\end{enumerate}
\end{corollary}
\begin{proof} Throughout, we assume $A=\emptyset$. It is immediate from Proposition \ref{P: family VA iff curve VA} that (2) implies (1) and (1) implies (3). To see that (3) implies (1), simply note that every plane curve in $X$ is a generic member of some $\emptyset$-definable family of plane curves, and the parameter space for such a family can be made stationary after passing to $\acl(\emptyset)$. So it is enough to show that (1) implies (2).

Now suppose (1) holds. As above, we can realize $C$ as a generic member of a $\emptyset$-definable family of plane curves, say $C_{t_0}\in\mathcal C=\{C_t:t\in T\}$. By compactness and the definability of dimension, we may assume that for each $t\in T$, for any distinct generics $x,y\in C_t$, we have $\dim(t/xy)\leq\dim T-2$ -- that is, that $\{C_t:t\in T\}$ is generically very ample. By Proposition \ref{P: family VA iff curve VA}, $\mathcal C$ becomes very ample after removing a non-generic $\emptyset$-definable subset of its graph. Since $t_0\in T$ is generic, only finitely many points of $C_{t_0}$ are removed, which is enough to prove (2).
\end{proof}

We conclude this section with a geometric characterisation of very ampleness: 

\begin{corollary}\label{very ample iff very ample family cor} Let $X$ be strongly minimal. Then $X$ is very ample if and only if there is a very ample definable family of plane curves in $X$.
\end{corollary}
\begin{proof}
If $X$ is very ample, by definition there is a plane curve in $X$ which is very ample over some set. Then by the previous corollary, there is a very ample family of plane curves in $X$.

If on the other hand there is a very ample family, then by the proposition there is a very ample plane curve over some set. Then by Corollary \ref{C: VA components} every strongly minimal compomnent of that plane curve is very ample over the same set, so that by definition $X$ is very ample.
\end{proof}

\subsection{Preservation Properties}

In this section we develop three basic preservation properties of very ampleness; we then pose a fourth (stronger) property as a question, and answer this question in the context of expansions of fields. To begin we point out, generalizing the preservation under naming parameters (Lemma \ref{L: Ind on parms}), that very ampleness is preserved under arbitrary expansions:

\begin{lemma} Any strongly minimal expansion of a very ample strongly minimal structure is very ample.
\end{lemma}
\begin{proof} Immediate by Corollary \ref{very ample iff very ample family cor}, since a very ample family is still definable in any expansion (technically we use here that the dimension of a definable set is unchanged in strongly minimal expansions, which is easy to check).
\end{proof}

Next we show that very ample plane curves are closed under `independent compositions,' in an appropriate sense. This will be used later on to show that very ample strongly minimal structures admit very ample families of plane curves of arbitrarily large dimension.

\begin{lemma}\label{L: composition} Let $X$ be strongly minimal and $A$-definable. Let $C$ and $D$ be plane curves in $X$, definable over $As$ and $At$, respectively. Assume that $s$ and $t$ are independent over $A$. If $C$ and $D$ are each very ample over $A$, then so is $D\circ C$.
\end{lemma}
\begin{proof}
Without loss of generality, assume $A=\emptyset$. Let $(x_1,z_1),(x_2,z_2)$ be distinct generics in $D\circ C$ over $st$. It will suffice to show that $\dim(st/x_1z_1x_2z_2)\leq\dim(st)-2$.

Without loss of generality, assume $x_1\neq x_2$. By definition of $D\circ C$ there are $y_1,y_2$ such that each $(x_i,y_i)\in C$ and each $(y_i,z_i)\in D$. Then each $(x_i,y_i)$ is generic in $C$ over $st$ (so also over $s$), which by the very ampleness of $C$ gives that $\dim(s/x_1y_1x_2y_2)\leq\dim(s)-2$.

Next, note that  $y_i$ and $z_i$ are interalgebraic over $t$ (for $i=1,2$), by the non-triviality of $D$ (see remark \ref{R: non trivial}). Thus $$\dim(st/x_1z_1x_2z_2)=\dim(sty_1y_2/x_1z_1x_2z_2)$$ $$=\dim(t/x_1z_1x_2z_2)+\dim(y_1y_2/tx_1z_1x_2z_2)+\dim(s/tx_1y_1z_1x_2y_2z_2).$$

In this last expression, note that the first term is at most $\dim(t)$, the second term is 0, and the third term is at most $\dim(s)-2$ as we have seen above. Thus, $$\dim(st/x_1z_1x_2z_2)\leq\dim(s)+\dim(t)-2,$$ which, by the independence of $s$ and $t$, is equivalent to the desired statement.
\end{proof}

Finally, we point out that very ampleness is preserved under finite-to-one functions:

\begin{lemma}\label{L: very ample images} Suppose $X$ and $Y$ are strongly minimal sets, and $f:X\rightarrow Y$ is a definable finite-to-one function. If $X$ is very ample, then so is $Y$.
\end{lemma}
\begin{proof}
Without loss of generality assume $X$, $Y$, and $f$ are $\emptyset$-definable. Let $C$ be a strongly minimal plane curve in $X$ which is very ample over some set $A$. Without loss of generality, we assume also that $A=\emptyset$. Let $c=\operatorname{Cb}(C)$. After editing finitely many points, we may assume $C$ is definable over $c$.

Now let $D$ be the image of $C$ in $Y^2$ (applying $f$ to each coordinate). So $D$ is definable over $c$. Since $f$ is finite-to-one, $D$ is moreover strongly minimal. We claim that $D$ is very ample in $Y$. To that end, suppose $z,w$ are distinct generics in $D$ over $c$. It will suffice to show that $z$ forks with $c$ over $w$. Now by assumption there are $x,y\in C$ with $f(x)=z$ and $f(y)=w$. Since $f$ is finite-to-one, each of the pairs $x,z$ and $y,w$ are interalgebraic. So by dimension considerations, $x$ and $y$ are generics in $C$. Moreover, since $z\neq w$, it follows that $x\neq y$. Thus, by the very ampleness of $C$, $x$ forks with $c$ over $y$. Then by interalgebraicity, $z$ forks with $c$ over $w$, as desired.
\end{proof}

 Lemma \ref{L: very ample images} does not extend to the preservation of  very ampleness  under \textit{non-orthogonality}, as can be seen from Example \ref{canonical example ex}. Indeed, the (non-very ample) structure $\mathcal M$ from Example \ref{canonical example ex} is non-orthogonal to the strongly minimal set $K/\sim$, where $\sim$ is the relation $x^2=y^2$ on $K$; but the induced structure on $K/\sim$ is a pure algebraically closed field, so very ample by Example \ref{acf line ex}. 

As we will see later, the main reason very ampleness does not pass from $K/\sim$ to $\mathcal M$ in Example \ref{canonical example ex} is that $K$ is not \textit{internal} to $K/\sim$ (in the sense of the reduct,  $\mathcal M$). We may thus ask the following:

\begin{question}\label{Q: VA under internality} Is very ampleness preserved under internality? That is, if $X$ is very ample, and $Y$ is strongly minimal and internal to $X$, must $Y$ be very ample?
\end{question}

To conclude this subsection, we now answer Question \ref{Q: VA under internality} whenever $X$ has a definable field structure. We will later improve on this result in Corollary \ref{C: field va int equivalence} and in the proof of Proposition \ref{P: going up}. In order to smoothly recall the setting in this later proposition, we present the result here as a corollary of two general statements (Lemma \ref{L: field e of i} and Proposition \ref{P: field va}). 

\begin{lemma}\label{L: field e of i} Every strongly minimal expansion of an algebraically closed field eliminates imaginaries.
\end{lemma}
\begin{remark} By \textit{expansion} above, we are assuming the field structure is part of the signature of $(K,+,\cdot,\dots)$. In general, for a strongly minimal structure, $\CM$, admitting  a definable field structure (with universe $M$), elimination of imaginaries will hold after naming any set of parameters defining a field structure on $M$.
\end{remark}
\begin{proof} Let $(K,+,\cdot,\dots)$ be an expansion of an algebraically closed field. Because of the field structure, it is automatic that $\acl(\emptyset)$ is infinite. So by weak elimination of imaginaries in strongly minimal structures, it suffices to code finite sets: that is, any such expansion $(K,+,\cdot,\dots)$ has elimination of imaginaries if and only if there are $\emptyset$-definable injections $K^{(n)}\rightarrow K^m$ for all $n$ (and some $m$ depending on $n$), where $K^{(n)}$ denotes the $n$-th symmetric power of $K$. But such functions already exist in the pure field $(K,+,\cdot)$, so we are done. 
\end{proof}

\begin{proposition}\label{P: field va} Let $(K,+,\cdot,\dots)$ be a strongly minimal expansion of an algebraically closed field, and let $X\subset K^n$ be a $\emptyset$-definable set of dimension $r\geq 1$. Let $\mathcal H=\{H_t:t\in T\}$ be the family of hyperplanes in $K^n$, and let $t\in T$ be generic. Then $\dim(H_t\cap X)=r-1$, and for any two distinct $x,y\in H_t\cap X$ generic over $t$, $x$ and $t$ fork over $y$.
\end{proposition}
\begin{proof} Throughout this proof, we freely use the fact that the strongly minimal expansion $(K,+,\cdot,\dots)$ preserves the dimensions of $(K,+,\cdot)$-definable sets, which is well known and easy to check. We also use that $T$ can be identified by a large subset of $\mathbb P^n(K)$ (the set of projective hyperplanes intersecting $K^n$); in particular, it follows that $T$ is stationary of dimension $n$. As is well-known, for any $x\in K^n$ the set $H^x$ (the hyperplanes through $x$) then has dimension $n-1$, and for $x\neq y$ the set $H^x\cap H^y$ has dimension $n-2$. 

Now we prove the proposition in two claims:
\begin{claim}
$\dim(H_t\cap X)=r-1$.
\end{claim} 
\begin{claimproof} Let $x\in X$ be generic, and let $s\in H^x$ be generic over $x$. Then $\dim(xs)=r+n-1$. Since $\dim(s)\leq n$, this forces $\dim(x/s)\geq r-1$. In particular $\dim(H_s\cap X)\geq r-1$.

Let us verify that  $\dim(H_s\cap X)=r-1$. Suppose not. Then there is $y\in H_s\cap X$ with $\dim(y/xs)=r$. Thus $\dim(xys)=2r+n-1$. But $\dim(xy)\leq 2r$, so $\dim(s/xy)\geq n-1$. As observed above, this is only possible if $x=y$; but clearly $x\neq y$, since $\dim(y/xs)=r\geq 1$.

So $\dim(H_s\cap X)=r-1$, and thus $\dim(x/s)\leq r-1$. Recalling that $\dim(xs)=r+n-1$, we conclude that $\dim(s)=n$ and $\dim(x/s)=r-1$. In particular, $s$ is generic in $T$. Since $T$ is stationary, $\tp(s)=\tp(t)$, so also $\dim(H_t\cap X)=r-1$. 
\end{claimproof}

\begin{claim} Let $x\neq y$ be generics in $H_t\cap X$ over $t$. Then $x$ forks with $t$ over $y$.
\end{claim}

\begin{claimproof}
By assumption $\dim(t)=n$, and by the previous claim $\dim(y/t)=r-1$. So $\dim(yt)=n+r-1$. But $\dim(t)\leq r$ since $y\in X$; and since $t\in H^y$, we get $\dim(t/y)\leq n-1$. So we must have $\dim(y)=r$ and $\dim(t/y)=n-1$. But since $t\in H^x\cap H^x$ we also have $\dim(t/xy)\leq n-2$, which proves the claim.
\end{claimproof}
\end{proof}

\begin{corollary}\label{C: field va} Let $K$ be a strongly minimal expansion of a field, and let $X$ be a strongly minimal set internal to $K$. Then $X$ is very ample.
\end{corollary}
\begin{proof} Adding parameters if necessary, we may assume $X$ is $\emptyset$-definable. By Lemma \ref{L: field e of i}, we may assume $X\subset K^n$ for some $n$. Then by Proposition \ref{P: field va} (applied to $X^2\subset K^{2n}$), there is a very ample plane curve in $X$. So $X$ is very ample.
\end{proof}

\subsection{When Very Ampleness can be Guaranteed}

 As we have seen, our first example of a non-very ample structure (Example \ref{canonical example ex}) still admits a very ample \textit{sort} -- that is, it becomes very ample after dividing by a definable equivalence relation. We thus think of $K$ of that example as very ample `up to a finite cover,' which is precisely the situation of Theorem B  of \cite{HrZil}. In an ideal world, this will hold for every non-locally modular strongly minimal structure:

\begin{question}\label{very ample sort question} Suppose $X$ is strongly minimal and not locally modular. 
\begin{enumerate}
    \item Can one always find a very ample strongly minimal set $Y$ internal to $X$? 
    \item More specifically, can one always find a definable equivalence relation $\sim$ on $X$ so that $X/\sim$ is strongly minimal and very ample?
\end{enumerate}
\end{question}

A conceivable approach to Question \ref{very ample sort question}(2) is as follows: first, show that there is a two-dimensional family $\{C_t:t\in T\}$ of plane curves in $\mathcal M$ such that the relation $I(x,y)$ on $M^2$, given by $I(x,y)\leftrightarrow|C^x\cap C^y|=\infty$, is in some sense induced by an equivalence relation $\sim$ on $M$; then, factor out this relation to form the sort $M/\sim$, and show that the image of $\{C_t:t\in T\}$ in $M/\sim$ (projecting each coordinate) is a very ample family. Later in this paper, we will see that this strategy is successful for strongly minimal groups (Theorem \ref{very ample group thm}). In the general case, it seems hard to show that the relation $I$ comes from an equivalence relation. However, we can show two related weaker statements. Proposition \ref{relation on M prop} shows that $I$ can indeed be taken to be induced by a relation on $M$ (rather than $M^2$):

\begin{proposition}\label{relation on M prop} Let $X$ be strongly minimal and not locally modular. Then there are a parameter set $A$, a 2-dimensional $A$-definable  faithful family $\{C_t:t\in T\}$ of plane curves in $X$, and a one-dimensional $A$-definable set $J\subset M^2$, such that for any two generics $x=(x_1,x_2),y=(y_1,y_2)\in M^2$ over $A$, if $C^x\cap C^y$ is infinite then either $J(x_1,y_1)$ or $J(x_2,y_2)$ holds.
\end{proposition}
Since the above proposition will not play a role in the sequel, we leave the proof as a simple exercise to the interested reader. \\

Proposition \ref{e of i prop} shows that, as one might expect from the above discussion, non-locally modular strongly minimal structures with elimination of imaginaries are very ample:

\begin{definition} Let $X$ be strongly minimal and $A$-definable. We say that $X$ \textit{eliminates imaginaries over $A$} if for any $a_1,\dots,a_n\in X$ and any tuple $c\in\dcl(A\bar a)$, there is a finite sequence $b_1,\dots,b_m\in X$ with $\dcl(A\bar b)=\dcl(Ac)$.
\end{definition}

\begin{proposition}\label{e of i prop} Let $X$ be strongly minimal, $A$-definable, and not locally modular. If $X$ eliminates imaginaries over $A$, then $X$ is very ample. 
\end{proposition}

\begin{proof} By non-local modularity, there are a set $B\supset A$ and a strongly minimal plane curve $C\subset X^2$, such that $\dim(c/B)=2$, where $c=\operatorname{Cb}(C)$. By stable embeddedness of $X$, there are $a_1,\dots,a_n\in X$ with $c\in\dcl(A\bar a)$. So by elimination of imaginaries, there are $b_1,\dots,b_m\in X$ with $\dcl(A\bar b)=\dcl(Ac)$, and in particular $\dcl(B\bar b)=\dcl(Bc)$. So $\dim(\bar b/B)=2$, and we may extract a two element basis for $\bar b$ over $B$. Without loss of generality, let us assume $\dim(b_1b_2/B)=2$. So $b_1b_2\in\dcl(Bc)$ and $c\in\acl(Bb_1b_2)$.

Now let $x\in M^2$ be generic in $C$ over $Bc$. So $\dim(cx/B)=3$. Since $c\notin\acl(B)$, $x$ forks with $c$ over $B$, which implies that $\dim(x/B)=2$, and thus $\dim(c/Bx)=1$. By interalgebraicity, $\dim(b_1b_2/Bx)=1$. Let $p=\operatorname{stp}(b_1b_2/Bx)$; so $p$ is a minimal type in $X^2$. Let $S$ be a strongly minimal plane curve whose generic type is $p$. So $(b_1,b_2)$ is generic in $S$ over $\acl(Bx)$.

\begin{claim} $S$ is very ample in $X$ over $B$.
\end{claim} 
\begin{claimproof} Let $(b_1',b_2')$ be another generic of $S$ over $\acl(Bx)$. So $(b_1',b_2')\models p$. Then there is some $c'$ with $\tp(b_1b_2cx/B)=\tp(b_1'b_2'c'x/B)$. Since $b_1b_2\in\dcl(Bc)$, and $(b_1',b_2')\neq(b_1,b_2)$, it follows that $c\neq c'$. Thus $c$ and $c'$ code nonparallel minimal types. Since $x$ realizes both we get $\dim(x/cc')=\dim(x/Bcc')=0$, so that $x$ and $c'$ fork over $Bc$.
\end{claimproof}

We have shown that $X$ admits a very ample strongly minimal plane curve,  
so $X$ is very ample.
\end{proof}
\begin{remark} The formulation of the above proof using types is convenient, but may obscure the geometric idea. Let us now explain the proof in more geometric terms: we first use non-local modularity to find a faithful two-dimensional family of plane curves; we then use elimination of imaginaries to reparametrize this family by a definable set $T\subset M^n$, and by taking a projection $M^n\rightarrow M^2$ we further reparametrize (up to losing generic strong minimality of the curves) by $M^2$. We then `dualize' the resulting family (interchanging points and curves), and observe that the `dual' notion of faithfulness is precisely very ampleness.
\end{remark}

\begin{remark} The converse of Proposition \ref{e of i prop} is false: Let $\mathcal M=(M,\dots)$ be the full $\mathbb C$-induced structure on a smooth irreducible curve $M$ over $\mathbb C$ of genus $g>0$. Then $\mathcal M$ is very ample by Corollary \ref{C: field va}, but does not eliminate the imaginary $\mathbb C$ over any set of parameters. Let us elaborate: first, by the Zilber trichotomy (\cite{HaSu} or \cite{CasACF0}), $\mathcal M$ interprets a set isomorphic to $\mathbb C$, and thus interprets a smooth curve $X$ of genus 0. If $X$ were embeddable into some $M^n$, then by composing with a projection one could obtain a non-constant rational map $X\rightarrow M$, which (by completing, normalizing, and then comparing the genus of each curve) is easily seen to contradict the Riemann-Hurwitz formula.
\end{remark}

\section{Applications of Very Ampleness}

We now turn toward applications of very ampleness, with an eye toward improved versions of the Zilber trichotomy for structures interpreted in algebraically closed fields. Throughout this section and Section \ref{S: applications} statements will often include the disclaimer "assuming Zilber's trichotomy" (for structures interpretable in algebraically closed fields). At the time of writing of this paper, the status of Zilber's trichotomy for such structures is as follows: 
\begin{enumerate}
    \item The second author and D. Sustretov have a proof of the trichotomy in the one dimensional case, \cite{HaSu}. 
    \item The first author has a proof of the trichotomy in characteristic 0 (for all dimensions), \cite{CasACF0}. This result is independent of \cite{HaSu}. 
    \item These results are yet to be published. 
    \item The first author and J. Ye have a proof of the trichotomy in positive characteristic (mostly, similar to the proof of \cite{CasACF0}). This proof is not yet publicly available. 
\end{enumerate}

\subsection{Very Ampleness and Internality}

The following proposition summarizes the effect of very ampleness in comparing strongly minimal sets. In fact, almost every application we make of very ampleness will, ultimately, follow from this result:

\begin{proposition}\label{internality prop} Let $X$ and $Y$ be strongly minimal sets, and assume $X$ is very ample. If $X$ is non-orthogonal to $Y$, then $X$ is internal to $Y$.
\end{proposition}

\begin{proof}
    By non-orthogonality, there is a strongly minimal set $C\subset X\times Y$ such that both projections $C\rightarrow X$ and $C\rightarrow Y$ are finite-to-one. Without loss of generality (by adding finitely many points to $C$) assume that $C\rightarrow X$ is surjective. We then view $C$ as a multivalued function $X\rightarrow Y$, mapping each $x\in X$ to the finite set $\{y\in Y:(x,y)\in C\}$. In particular, it follows that one can find a sort $S$ in $Y^{\textrm{eq}}$, and a finite-to-one definable map $f:X\rightarrow S$. Absorbing parameters into the language, we may assume that $f$ is $\0$-definable. 
    
    Now by assumption there is a plane curve $C\sub X^2$ which is very ample over some set $A$. Applying $f$ coordinate-wise, we obtain a finite-to-one function $\tilde f: C\to S^2$. Note that since $f$ is $\0$-definable, if $\tilde f(x_1,x_2)=\tilde f(y_1,y_2)$ then $(x_1,x_2),(y_1,y_2)$ are interalgebraic over $\0$. So by condition (3) of Lemma \ref{L: Ind on parms}, there are no distinct generics $(x_1,x_2), (y_1,y_2)\in C$ over $A$ such that $\tilde f(x_1,x_2)=\tilde f (y_1,y_2)$. In other words, $\tilde f$ is generically injective. Letting $T=\im(\tilde f)$ (potentially minus a finite set) and inverting $\tilde f$, we can then extract a definable function $g: T\to C$ with cofinite image. Then composing with a coordinate projection $\pi$, the resulting map $\pi\circ g: T\to X$ has cofinite image in $X$. This shows that $X$ is internal to $T$, and since $T\subset S^2$, $T$ is internal to $Y$, thus $X$ is moreover internal to $Y$.
\end{proof}

Before moving on we give two quick applications. First, the following is immediate from Proposition \ref{internality prop}, Lemma \ref{L: field e of i}, and Corollary \ref{C: field va}:

\begin{corollary}\label{C: field va int equivalence} Let $K$ be a strongly minimal expansion of an algebraically closed field, and let $X$ be a strongly minimal set non-orthogonal to $K$. Then the following are equivalent:
\begin{enumerate}
    \item $X$ is very ample.
    \item $X$ is internal to $K$.
    \item There is a definable embedding of $X$ into some $K^n$.
\end{enumerate}
\end{corollary}

It follows from this last corollary, that if $X$ is very ample and interprets an algebraically  closed field, $K$, then any strongly minimal set $Y$ interpretable in $X$ is very ample.  Indeed,  by the corollary $X$ is internal to $K$, and since $Y$ is internal to $X$, $Y$ too is internal to $K$, so by the corollary,  again, $Y$ is very ample. Thus, in particular, a counter-example to Question \ref{Q: VA under internality} must be a counter-example to Zilber's trichotomy. \\

Next, we give a partial answer to Question \ref{very ample sort question} for strongly minimal sets non-orthogonal to very ample sets:

\begin{corollary}\label{C: partial va sort}
Let $X$ and $Y$ be strongly minimal and non-orthogonal, and assume $Y$ is very ample. Then:
\begin{enumerate}
    \item There is a very ample strongly minimal set internal to $X$.
    \item If, moreover,  $Y$  eliminates imaginaries (over some parameter set), then there is a definable equivalence relation $\sim$ on $X$ with finite classes such that $X/\sim$ is very ample.
\end{enumerate}
\end{corollary}
\begin{remark} In particular, it follows from Corollary \ref{C: partial va sort} that every non-locally modular strongly minimal structure satisfying the Zilber trichotomy admits a positive answer to both parts of Question \ref{very ample sort question}.
\end{remark}
\begin{proof} (1) is immediate, because by Proposition \ref{internality prop} $Y$ is internal to $X$. So we need only prove (2). For that, assume that $Y$ eliminates imaginaries over some set. By non-orthogonality, one can find a finite-to-one map $f:X\rightarrow S$, where $S$ is a sort in $Y^{\textrm{eq}}$. By elimination of imaginaries we may assume $S\subset Y^n$ for some $n$, so we have $f:X\rightarrow Y^n$. Then composing with an appropriate projection $Y^n\rightarrow Y$ we obtain a finite-to-one map $X\rightarrow Y$, which shows that (up to editing finitely many points) $Y$ is a quotient of $X$.
\end{proof}

\subsection{Very ample Zariski Geometries}\label{Ss: ZG}
In the present section we show that our definition of very ampleness is equivalent, in the context of Zariski Geometries, to the original definition of Hrushovski and Zilber. The results of this section are not needed in the sequel, and the reader only interested in the application may safely skip to the next section.  

Recall that a strongly minimal Zariski Geometry\footnote{We are referring to the term as appearing in \cite{HrZil} and not to the more general one of \cite{ZilZariski}.} is a structure $\CZ$, in a language equipping it  with a compatible system of Noetherian topologies on $Z^n$ (all $n$), such that $\CZ$ has quantifier elimination and satisfies the dimension theorem (see \cite{HrZil} for the details). The Zariski Geometry $\CZ$ is ample if it is not locally modular; it is then shown \cite{HrZil}[Theorem B] that any ample Zariski geometry interprets an algebraically closed field $K$. To further identify $\CZ$ with a smooth curve over $K$, the notion of very ampleness is introduced. We now give the definition (paraphrased) as it appears in \cite{HrZil}. For the sake of clarity, we will refer to this version as Z-very ample.

\begin{definition} Let $\CZ$ be a Zariski geometry with universe $Z$. Then $\CZ$ is \textit{Z-very ample} if there are an irreducible closed set $E\subset Z^n$ for some $n$, and an irreducible closed set $C\subset E\times Z^2$, such that the following hold:
\begin{enumerate}
\item For generic $e\in E$, the set $C(e)$ is irreducible and one-dimensional.
\item For any distinct $a,b\in Z^2$ there is some $e\in E$ such that $C(e)$ contains exactly one of $a,b$.
\end{enumerate}
\end{definition}

The main theorem of \cite{HrZil}, equivalently stated, then asserts that every Z-very ample Zariski geometry is isomorphic (as a Zariski geometry) to a smooth curve over an algebraically closed field.

It is not directly mentioned in \cite{HrZil} that the converse holds -- namely, that every smooth curve is very ample. We thank E. Hrushovski for outlining the following proof.

\begin{lemma}\label{L: curves are zva}
    Let $C$ be a smooth algebraic curve over an algebraically closed field. Then $C$ is $Z$-very ample. 
\end{lemma}
\begin{proof} Since $C$ is a curve it is quasi-projective; so we may assume $C\subset\mathbb P^n(K)$ for some $n$. Our goal is to find a family $C\subset E\times Z^2$ consisting of all hyperplane intersections in $Z^2$ (that is, sets $H_t\cap Z^2$ where $H_t$ is a hyperplane in $\mathbb P^n(K)$). In this case (1) follows by Bertini's Theorem, and (2) follows since any two points in projective space can be separated by a hyperplane.

The only difficulty here is finding such a family which is irreducible, closed, and parametrized by a closed irreducible set in a power of $Z$ (rather than the natural parametrization by $\mathbb P^n(K)$). For this, it will suffice to find a surjective morphism $Z^n\rightarrow\mathbb P^n(K)$, since then we can lift the natural parametrization from $\mathbb P^n(K)$ to $Z^n$.

To find such a morphism, we note that by composing with a surjection $(\mathbb P^1)^n\rightarrow\mathbb P^n$ (which exists because the former is a projective variety of dimension $n$), it will moreover suffice to find a surjection $C\rightarrow\mathbb P^1(K)$.

Finally, let us build a surjection $C\rightarrow\mathbb P^1(K)$. First choose any (dominant) morphism $f:C\rightarrow\mathbb P^1(K)$, and suppose $f(C)$ misses $m$ points of $\mathbb P^1(K)$, say $x_1,\dots ,x_m$. Then the map $g\circ f:C\rightarrow\mathbb P^1(K)$ is surjective, where $g:\mathbb P^1(K)\rightarrow\mathbb P^1(K)$ is any ramified cover with generically $m+1$-sized fibers that does not ramify at any of $x_1,\dots,x_m$.
\end{proof}

It follows by \cite{HrZil} and Lemma \ref{L: curves are zva} that a Zariski geometry $\CZ$ is Z-very ample if and only if it is isomorphic to a smooth curve over an algebraically closed field. Our goal now is to show that the same statement holds using our notion of very ampleness in place of Z-very ampleness. To do this, we first need to inspect the proof in \cite{HrZil}, in order to extract the precise use of Z-very ampleness and replace it with a use of very ampleness. In fact, the proof in \cite{HrZil} is achieved by separately proving the following three facts:

\begin{fact}\label{F: ZG proof overview} Let $\CZ$ be a Zariski geometry with universe $Z$.
\begin{enumerate}
    \item If $\CZ$ is not locally modular then $\CZ$ interprets an algebraically closed field (p. 47, first paragraph of the proof of Theorem A).
    \item If $\CZ$ interprets the algebraically closed field $K$, and $\CZ$ is Z-very ample, then $Z$ is internal to $K$ (from (1) to the end of the paragraph splitting p. 47-48).
    \item If $\CZ$ interprets the algebraically closed field $K$, and $\CZ$ is internal to $K$, then $\CZ$ is isomorphic to a smooth curve over $K$ (p. 48, the rest of the proof of Theorem A).
\end{enumerate}
\end{fact}

We now show:

\begin{proposition}\label{P: zg va equiv}
    Let $\mathcal Z$ be a Zariski Geometry. Then the following are equivalent: 
    \begin{enumerate}
        \item $\CZ$ is very ample.
        \item $\mathcal Z$ is $Z$-very ample. 
        \item $\mathcal Z$ is isomorphic to a smooth curve over an algebraically closed field.
    \end{enumerate}
\end{proposition}
\begin{proof}
    As stated above, the equivalence of (2) and (3) is \cite{HrZil}[Theorem A] and Lemma \ref{L: curves are zva}. Moreover, the implication (3) $\Rightarrow$ (1) follows from Corollary \ref{C: field va}. So it will suffice to show (1) $\Rightarrow$ (3).

    Assume (1). Then by Lemma \ref{L: very ample implies nlm}, $\CZ$ is not locally modular. So by Fact \ref{F: ZG proof overview}(1), $\CZ$ interprets an algebraically closed field $K$. Then by Proposition \ref{internality prop}, $Z$ is internal to $K$. So by Fact \ref{F: ZG proof overview}(3), $\CZ$ is isomorphic to a smooth curve over $K$. Thus (3) holds, and we are done.
\end{proof}

\subsection{Very Ample Strongly Minimal ACF-Relics}

We now turn toward applications of very ampleness to structures interpretable in algebraically closed fields. In \cite[Theorem 7.1]{Lov} Loveys gives a complete list (up to definable finite covers) of all locally modular non-trivial reducts of the algebraically closed field $\Cc$. If we are  interested in the classification of non-locally modular reducts of $\Cc$ only up to finite covers, then Rabinovich' theorem, \cite{Ra}, is the final word. But one could hope for a more precise classification of such reducts. One natural question is the identification of reducts that are not proper. To state this problem in somewhat greater generality, it is convenient to have: 
 \begin{definition}
     Let $\CN$ be any structure. An \emph{$\CN$-relic} is a reduct, $\mathcal X$, of the structure induced on some $\CN$-definable set $X$. The relic $\mathcal X$ is \emph{full} in $\CN$ if any $\CN$-definable subsets of $X^m$ (any $m$) is $\mathcal X$-definable. 
 \end{definition}
 
In a recent unpublished work, the first author and M. Tran give a fairly general condition for a reduct $\CM$ of $\Cc$ whose atomic sets are polynomial functions to be full, though a complete identification of all full reducts of $\Cc$ is not available. In the next two subsections we show, roughly, that the fullness of an ACF-relic is equivalent to very ampleness. In fact, in the main result of the current subsection, we show an equivalence between several properties of a strongly minimal ACF-relic, notably very ampleness, fullness, and bi-interpretability with the field.

As a warm-up, we begin by pointing out two straightforward corollaries of Proposition \ref{internality prop} in the context of algebraically closed fields, which highlight the enhancements very ampleness brings to the statement of the Zilber trichotomy. First, in analogy with the main theorem of \cite{HrZil} on Zariski Geometries, we show that (assuming the Zilber trichotomy) every very ample strongly minimal ACF-relic is isomorphic (outside a finite set) to a curve (rather than being a finite cover of a curve):

\begin{corollary}\label{non-orthogonal implies internal cor} Suppose $X$ is strongly minimal and non-orthogonal to a pure algebraically closed field $K$. If $X$ is very ample, then after removing finitely many points $X$ is isomorphic to an irreducible algebraic curve over $K$, with its full $K$-induced structure.
\end{corollary}
\begin{proof} By Corollary  \ref{C: field va int equivalence}, we may assume $X\subset K^n$ for some $n$. By the purity of $K$, the induced structure on $X$ is just the induced structure from the field. Finally, as a strongly minimal subset of $K^n$, $X$ is the union of an irreducible algebraic curve and a finite set. Thus, outside finitely many points, we identify $X$ with an irreducible curve, as desired.
\end{proof}

Next, we show (again assuming the Zilber trichotomy) that all very ample strongly minimal ACF-relics are full:

\begin{corollary}\label{C: va implies full} Let $K$ be an algebraically closed field, and let $\mathcal M=(M,\dots)$ be a strongly minimal $K$-relic which satisfies Zilber's trichotomy. If $\mathcal M$ is very ample then $\mathcal M$ is full in $K$. 
\end{corollary}
\begin{proof}
By Lemma \ref{L: very ample implies nlm}, $\mathcal M$ is not locally modular. So by the Zilber Trichotomy, there is a field $F$ isomorphic to $K$ which is interpreted in $\mathcal M$. By \cite{PoiFields}, there is a $K$-definable isomorphism $K\leftrightarrow F$. In particular, the induced structure on $F$ from $K$ is that of a pure field. Since $\mathcal M$ defines the field operations on $F$, it follows that the induced structure on $F$ from $\mathcal M$ is also a pure field.

Now, since $\mathcal M$ interprets $F$, $M$ and $F$ are non-orthogonal in $\mathcal M$. So by Corollary \ref{C: field va int equivalence}, there is an $\mathcal M$-definable bijection $f:M\rightarrow X$, where $X\subset F^n$ is some $F$-constructible set. 

Now, to prove the corollary, let $Y\subset M^k$ be $K$-constructible. Then $f(Y)\subset X^k$ is $K$-definable, and thus $F$-constructible by the purity of $F$ in $K$. Since $\mathcal M$ interprets $F$, this implies that $f(Y)$ is $\mathcal M$-definable, and thus (inverting $f$) so is $Y$. Thus $\mathcal M$ defines all constructible subsets of powers of $M$.
\end{proof}

To summarize, we give the following characterization of very ampleness for strongly minimal ACF-relics:

\begin{theorem}\label{very ample acf char thm} Let $K$ be an algebraically closed field, and let $\mathcal M=(M,\dots)$ be a non-locally modular strongly minimal $K$-relic which satisfies Zilber's trichotomy. Then the following are equivalent:
\begin{enumerate}
    \item $\mathcal M$ is very ample.
     \item $\mathcal M$ is full in $K$.
    \item $\mathcal M$ is internal to every infinite field it interprets.
    \item $\mathcal M$ is internal to some infinite field it interprets.
    \item $\mathcal M$ is bi-interpretable with $K$ over parameters.
    \item Up to deleting a finite set, $\mathcal M$ is isomorphic to an irreducible algebraic curve over $K$, with its induced structure from $K$.
\end{enumerate}
\end{theorem}
\begin{remark} Note, in particular, that fullness a priori depends on the interpretation of $\mathcal M$ in $K$, while the other five conditions do not. So it follows that fullness is a property only of the abstract structure $\mathcal M$.
\end{remark}

\begin{proof} The equivalence of (1) and (2) is given by Corollaries \ref{C: va implies full} and \ref{C: field va}. We show that each of the other items is equivalent to (1).
\begin{itemize}
    \item[$(1)\Rightarrow(3)$] Immediate by Proposition \ref{internality prop}.
    \item[$(3)\Rightarrow(4)$] By the Zilber trichotomy $\mathcal M$ interprets an infinite field, which is enough.
    \item[$(4)\Rightarrow(5)$] Note that any infinite field interpreted in an algebraically closed field $K$ is $K$-definably isomorphic to $K$ by the main result of \cite{PoiFields}. So it suffices to find a field $F$ interpreted in $\mathcal M$ (and thus isomorphic to $K$), and a definable set $X\subset F^n$ for some $n$, such that the induced structure on $X$ from $F$ is $\mathcal M$-definably isomorphic to $M$. For this last condition, it is sufficient to show that $M$ and $X$ are in $\mathcal M$-definable bijection, and that the induced structures on $X$ from $F$ and $\mathcal M$ coincide.
    
    Now we can satisfy the above requirements using (4). That is, by (4) there are a field $F$ interpretable in $\mathcal M$, and an $\mathcal M$-definable bijection between $M$ and an $\mathcal M$-definable set $X$ in $F^{\textrm{eq}}$. As above, $K$ and $F$ are $K$-definably isomorphic, which shows that the induced structure on $F$ from $K$ (and thus also from $\mathcal M$) is a pure field isomorphic to $K$. In particular, we can assume $X$ is constructible over $F$, and that the induced structure on $X$ from $\mathcal M$ is the same as the induced structure from $F$. 
    \item[$(5)\Rightarrow(6)$] Given a bi-interpretation, and arguing as in the previous point, we obtain an $\mathcal M$-definable field $F$ and an $\mathcal M$-definable bijection between $M$ and some $F$-constructible set $X$. As above, the induced structure on $F$ from both $\mathcal M$ and $K$ is a pure field isomorphic to $K$. Thus, the image of any $\mathcal M$-definable set $Y\subset M^k$ is an $F$-constructible subset of $X^k$; moreover, the preimage of any $F$-constructible subset of $X^k$ is $\mathcal M$-definable, so is an $\mathcal M$-definable subset of $M^k$. We have thus shown that the induced structures on $M$ and $X$ from $\mathcal M$ are isomorphic. So it suffices to replace $\mathcal M$ with $X$ in proving (6). But $X$ is then a strongly minimal set in the pure field $F$ with the full induced structure from $F$, so is thus (up to deleting a finite set) isomorphic to the full induced structure on an irreducible curve over $F$. Finally, since $F$ and $K$ are isomorphic, $X$ is isomorphic to an analogous structure over $K$, as desired.
    \item[$(6)\Rightarrow(1)$] If (6) holds, then $\mathcal M$ is isomorphic to a full relic, and we can apply Corollary \ref{C: field va}.
    \end{itemize}
\end{proof}

\subsection{More on full ACF-relics}
It follows from Theorem \ref{very ample acf char thm} that strongly minimal relics of algebraically closed fields are full precisely when they are very ample. In this section we give a similar characterization of fullness for a general (non-strongly minimal) relic. In this case the argument will be a little more delicate.

\begin{theorem}\label{T: higher relics}
    Let $K$ be an algebraically closed field, and assume the Zilber trichotomy holds for strongly minimal $K$-relics. Let $\CM=(M,\dots)$ be a $K$-relic. Then the following are equivalent:
    \begin{enumerate}
        \item $\CM$ is full.
        \item $\CM$ is almost strongly minimal, and every strongly minimal set in $\CM$ is very ample.
    \end{enumerate}
\end{theorem}

\begin{proof}
    First, assume $\mathcal M$ is full; we show (2). Let $X$ be any strongly minimal set in $\mathcal M$. Then $X$ is full (as a relic with its induced structure from $\mathcal M$), so very ample by Corollary \ref{C: field va}. Now applying the Zilber trichotomy to any such $X$ (\cite{HaSu} -- noting that by fullness we have $\dim_K(X)=1$ and by Lemma \ref{L: very ample implies nlm} $X$ is not locally modular), it follows that $\mathcal M$ interprets a copy $K'$ of $K$. It then follows (potentially after adding a small set of parameters) that $M\sub \acl_{\mathcal M}(K')$. Indeed, in $K$ there is a definable bijection between $K'$ and $K$, so $\acl_K(K')\supseteq M$, and we conclude by the fullness of $\CM$. 
    
    Now assume $\mathcal M$ is almost strongly minimal, and all strongly minimal sets in $\mathcal M$ are very ample. Then, by Corollary \ref{C: va implies full} and the Zilber trichotomy for $K$-relics, $\CM$ interprets a (strongly minimal) algebraically closed field $F$. After adding parameters to both $\mathcal M$ and $K$, we will assume that $F$ is $\0$-definable in both structures. We will also assume that $\mathcal M$ (and all of its basic relations) are $\emptyset$-definable in $K$ (this may require passing to an elementary extension to preserve saturation). 
    
    Now it follows by almost strong minimality that every strongly minimal set in $\mathcal M$ is non-orthogonal to $F$; It then moreover follows by Proposition \ref{internality prop} that every strongly minimal set in $\mathcal M$ is internal to $F$. Our plan is to show more generally that every $\CM$-definable set is internal to $F$, using an induction on dimension. Throughout, we use the following two lemmas:
    
    \begin{lemma}\label{L: internal implies full} Let $X$ be any $\mathcal M$-definable set, and let $\mathcal X$ denote the $\CM$-induced structure on $X$  (so $\mathcal X$ is naturally a $K$-relic). If $X$ is internal to $F$, then $\mathcal X$ is full in $K$.
    \end{lemma}
    \begin{proof} The idea is well known and has already been seen in Corollary \ref{C: va implies full} and Theorem \ref{very ample acf char thm}. First, note that since there is a $K$-definable isomorphism $K\rightarrow F$, the induced structure on $F$ from $K$ is a pure field. Now, by internality (and the previous sentence), there is an $\mathcal M$-definable bijection $g: X\to  Y$ for some $F$-constructible set $Y$. Let $Z\subset X^k$ be any $K$-definable set. Then $g(Z)\subset Y^k$ is $F$-constructible, so $\mathcal M$-definable. And since $g$ is $\mathcal M$-definable, it follows that $g^{-1}(g(Z))=Z$ is also $\mathcal M$-definable. Finally, by definition of $\mathcal X$ this shows that $Z$ is  $\mathcal X$-definable, proving that $\mathcal X$ is full.
    \end{proof}

    \begin{lemma}\label{L: rank preserving} Let $X$ be any $\mathcal M$-definable set. Then $\dim_{\mathcal M}(X)=\dim_K(X)$.
    \end{lemma}
    \begin{proof} By almost strong minimality, there is an $\mathcal M$-definable finite-to-one map $f:X\rightarrow F^k$ for some $k$. So $X$ is in $\mathcal M$-definable finite correspondence with an $F$-constructible set; since both $\dim_{\mathcal M}$ and $\dim_K$ respect finite correspondences, we may assume $X$ is itself $F$-constructible. In this case both $\dim_{\mathcal M}(X)$ and $\dim_K(X)$ are equal to $\dim_F(X)$, because (arguing as in Lemma \ref{L: internal implies full}) the induced structures on $F$ from $F$, $K$, and $\mathcal M$ all coincide.
    \end{proof}

    \begin{nota} In light of Lemma \ref{L: rank preserving}, for the rest of the proof we will drop all subscripts when denoting dimensions of definable sets.
    \end{nota}
    
    We now proceed with the inductive step (toward showing that every $\CM$-definable set is $F$-internal). So let $X$ be any $\mathcal M$-definable set, say of dimension $r$, and assume that every $\mathcal M$-definable set of smaller dimension is internal to $F$. Since a finite union of sets internal to $F$ is clearly internal to $F$, we may assume $X$ is stationary. Since we have already observed that strongly minimal sets are internal to $F$, we may then assume that $r\geq 2$. 
    
    By almost strong minimality, we get an $\CM$-definable finite-to-one map $f: X\to F^n$ for some $n$. Since $X$ is stationary, such a map must be generically $k$-1 for some $k$. We fix onece and for all such an $f$ for which $k$ is minimal. 

    \begin{nota} In the rest of this proof, we will use canonical bases in the senses of both $\mathcal M$ and $K$,  denoted  $\operatorname{Cb}_{\mathcal M}$ and $\operatorname{Cb}_K$, respectively.
    \end{nota}
    
    \begin{claim}
    Let $D$ be a stationary subset of $X$ of dimension $r-1$, and assume that $\mathrm{Cb}_{\mathcal M}(D)\notin\acl_{\mathcal M}(\0)$. Then $D$ is almost equal to a union of $f$ fibres.
    \end{claim}
    \begin{claimproof}
        Since $\operatorname{Cb}_{\mathcal M}(D)\notin\acl_{\mathcal M}(\0)$, we can express $D$ (up to an error of smaller rank)  as a generic member $D_{t_0}$ of a rank $1$ family of rank $r-1$ subsets of $X$, say $\{D_t\}_{t\in T}$ (potentially defined over additional parameters). By induction $T$ is $F$-internal. For generic $x\in  X$ consider the set $D^x:=\{t\in T: x\in D_t\}\subseteq T$. Since $T$ is $F$-internal, the set $D^x$ is coded in $F$, and we can append to the  map $f$ the data $x\mapsto \ulcorner D^x \urcorner$. Since we have chosen $f$ with minimal fibres, appending this information to $f$ cannot reduce the size of the fibres, implying that  if $x,y\in X$ are generic and $f(x)=f(y)$ then $D^x=D^y$. Now let $x$ be generic in $D=D_{t_0}$ over $t_0$, and let $y\in f^{-1}(f(x))$. An easy computation shows that $x$ and $y$ are generics in $X$, so from above $D^x=D^y$. By definition $t_0\in D^x$, so also $t_0\in D^y$, and thus $y\in D$. This proves the claim. 
    \end{claimproof}
    
    \begin{remark} Let $D$ be as in the claim. Then by the conclusion of the claim,  $f^{-1}f(D)$ almost coincides with $D$ (because $D\subset f^{-1}(f(D))$ always holds, and the claim shows that $f^{-1}f(x)\sub D$ for all $x\in D$ generic). 
    \end{remark}
    
    We now use the inductive hypothesis to show that the above claim holds even in the sense of $K$: 
    
    \begin{claim}
    Let $C\sub X$ be a constructible $K$-stationary set of dimension $r-1$ such that $\operatorname{Cb}_K(C)\notin\acl_K(\emptyset)$. Then $C$ almost coincides with a union of $f$-fibers.
    \end{claim}
    \begin{claimproof}
        Consider $f(C)\sub F^n$. Since $F$ is a pure algebraically closed field interpretable in $\CM$ it is full, thus $f(C)$ is $\CM$-definable and so is $f^{-1}f(C)$. Since $f$ is finite-to-one, $\dim(f^{-1}(C))=\dim(f(C))=r-1$. Then since $C$ is $K$-stationary, so is $f(C)$; clearly this implies that $f(C)$ is also $\mathcal M$-stationary.
        
        Let $D\sub f^{-1}f(C)$ be an $\CM$-stationary subset of dimension $r-1$. Clearly, then, $f(D)$ almost coincides with $f(C)$. If we can show that $\operatorname{Cb}_{\mathcal M}(D)\notin\acl_{\mathcal M}(\0)$, then by the previous claim $D$ almost coincides with $f^{-1}f(D)$, and thus also with $f^{-1}(f(C))$. So $C$ is generically contained in $D$ (because $C$ is automatically contained in $f^{-1}(f(C))$). Thus $\dim(C\cap D)=r-1$. But by the inductive hypothesis $D$ is $F$-internal, so by Lemma \ref{L: internal implies full} $D$ is full. Thus $C\cap D$ is an $\mathcal M$-definable $r-1$-dimensional subset of $D$; so by the stationarity of $D$, $C\cap D$ almost coincides with $D$. In other words, this shows that $C$ almost contains $D$; then combined with the above fact that $C$ is almost contained in $D$, we conclude that $C$ almost coincides with $D$. Then $C$ almost coincides with $f^{-1}(f(C))$ (because $D$ does), proving the claim in this case.
        
        So it remains to verify that $\operatorname{Cb}_{\mathcal M}(D)\notin\acl_{\mathcal M}(\0)$. But, as we have just seen, $f(C)$ and $f(D)$ are almost equal $F$-constructible sets which are stationary in both $\mathcal M$ and $K$. Since $F$ is a pure field in both $\mathcal M$ and $K$, it follows that $\operatorname{Cb}_{\mathcal M}(f(C))=\operatorname{Cb}_{\mathcal M}(f(D))=\operatorname{Cb}_K(f(C))=\operatorname{Cb}_K(f(D))$. Now, clearly, $\operatorname{Cb}_{\mathcal M}(D)$ and $\operatorname{Cb}_{\mathcal M}(f(D))$ are $\mathcal M$-interalgebraic; so if the former is in $\acl_{\mathcal M}(\0)$, then so is the latter. It clearly follows that $\operatorname{Cb}_{\mathcal M}(f(D))\in\acl_K(\0)$, and thus by above $\operatorname{Cb}_K(f(C))\in\acl_K(\0)$. But (as above with $D$ and $\mathcal M$) clearly $\operatorname{Cb}_K(C)$ and $\operatorname{Cb}_K(f(C))$ are interalgebraic, so this gives $\operatorname{Cb}_K(C)\in\acl_K(\0)$, contradicting the assumption in the claim. 
    \end{claimproof}

    We now aim to show that $X$ is internal to $F$:

    \begin{claim} $f$ is generically injective.
    \end{claim}
    \begin{claimproof} Throughout the proof of the claim, we work in the structure $K$ (so $\dim(a/A)$ means $\dim_K(a/A)$, and genericity means in the sense of $K$). By elimination of imaginaries in $K$, we may assume that $X\subset K^m$ for some $m$. We may also assume that $X$ and $f$ are $\0$-definable in $K$. We now apply Proposition \ref{P: field va}. Namely, let $\{H_s:s\in S\}$ be the family of hyperplanes in $K^m$, and let $s\in S$ be generic. Then $\dim(H_s\cap X)=r-1$, and for any distinct  $x,y\in H_s\cap X$ generic over $s$, $x$ forks with $s$ over $y$.

    Let $C\subset H_s\cap X$ be a constructible $K$-stationary set of dimension $r-1$. We may assume $C$ is $K$-definable over $\acl_K(s)$. Let $x\in C$ be generic over $\acl_K(s)$. Then $\dim(xs)=m+r-1$, while $\dim(x)\leq r$ and $\dim(s/x)\leq m-1$; so clearly $\dim(x)=r$ and $\dim(s/x)=m-1$. Note, then, that $\operatorname{Cb}_K(C)\notin\acl_K(\0)$, since otherwise we would have $\dim(x)=r-1$.

    Now, toward a proof of the claim, let $y\in X$ with $f(x)=f(y)$. It will suffice to show that $y=x$. Suppose not. By the previous claim, $C$ is almost equal to a union of $f$ fibres. So since $x\in C$ is generic, the fibre $f^{-1}(f(x))$ should be contained in $C$. In particular, we get $y\in C$, and thus $y\in H_s\cap X$. But since $f(x)=f(y)$, clearly $x$ and $y$ are interalgebraic over $\0$, so $y$ is also generic in $H_s\cap X$ over $s$. We conclude that $x$ must fork with $s$ over $y$; but this is impossible, because $x$ is algebraic over $y$, so $\tp(x/y)$ has no forking extensions.
    \end{claimproof}
    
    Now, by the most recent claim (and Lemma \ref{L: rank preserving}), there is an $\mathcal M$-definable large subset $X'\subset X$ such that $f$ is injective on $X'$. Then $X'$ is internal to $F$, witnessed by $f$. Since $X'$ is large in $X$ we get $\dim(X-X')<r$, so by the inductive hypothesis $X-X'$ is also internal to $F$. Then $X=X'\cup(X-X')$ is also internal to $F$, as desired.

    This completes the proof that every $\mathcal M$-definable set is internal to $F$. Now, in particular, this implies that the universe $M$ of $\mathcal M$ is internal to $F$. Then, by Lemma \ref{L: internal implies full} applied to $M$, we conclude that $\mathcal M$ is full. This completes the proof of the theorem.    
\end{proof}

\begin{remark} The reader may be tempted to replace `almost strongly minimal' with `uncountably categorical' in the statement of Theorem \ref{T: higher relics}. Let us point out, then, that almost strong minimality and unidimensionality are not equivalent for ACF-relics. Indeed, for any algebraically closed field, $K$, consider the two sorted relic $\CM$ defined as follows: we define the two sorts to be $M:=K^2$ and $U:=K$. We then equip $\mathcal M$ with the binary relation $V$ on $M$ interpreted as $\{(x,y), (x,z)): x,y,z\in K\}$; the map $f((x,y),(x,z))=y-z$ from $V$ to $U$; and the full field structure on $U$. (Essentially, we have defined a $K$-indexed family of copies of $K$, each of which has `forgotten where 0 is' but remembers an isomorphism with $K$ after fixing any element to be 0). Clearly, $\mathcal M$ is a $K$-relic. We leave it as an exercise to the interested reader to verify that $\mathcal M$ is uncountably categorical but not almost strongly minimal. 
\end{remark}

\section{Technical Results on Strongly Minimal Structures}

Here we prove two useful facts about general non-locally modular strongly minimal structures, without assuming very ampleness or Zilber's trichotomy. Though these are technical results needed for the sequel, we hope they will be seen as interesting in their own right. \textbf{Throughout this section $\CM$ is assumed to be a non-locally modular strongly minimal structure.} We assume throughout that $\acl(\0)$ is infinite, and that for each $k$ one can find a $\0$-definable, faithful family of plane curves $\{C_t:t\in T\}$ in $M$ where $\dim T=k$. Note that these conditions can be arranged by adding constants to the language. We point out, then, that the results of this section (Propositions \ref{P: large families} and \ref{P: transverse families}) are unchanged after dropping any added constants (except for the loss of $\0$-definability of the sets constructed in Proposition \ref{P: large families}). So if we do not insist on $\0$-definability, we can use these propositions freely in arbitrary non-locally modular $\mathcal M$. 

\subsection{Large Families Inside Any Definable Set}

It is well known that non-local modularity is characterized by the presence of arbitrarily large families of 1-dimensional subsets of $M^2$ (see also Conditions (A), (B) and (C) on pp. 85-87 of \cite{HrUnidim}). It is natural to wonder, then, whether there is anything special about 1-dimensional subsets of $M^2$: namely, suppose we are given any definable set $Z$ of dimension $n$; can we find arbitrarily large families of $m$-dimensional subsets of $Z$, for every $0<m<n$? 

On the face of it, the above-mentioned Condition (B1) of \cite{HrUnidim} seems to give a positive answer to our question. 
However, using the same measure of largeness used with plane curves (the dimension of the parameter space in a faithful parametrization) does not seem to really capture the right notion if $n>2$. For example, given arbitrarily large families of plane curves, one can trivially find high-dimensional faithful families of curves in $M^3$ by simply embedding $M^2$ into $M^3$. Clearly, though, it would be nice to have a construction which works on the level of types, so that (for example) a generic element of a generic member of the family is generic in $M^3$ (or, more generally, in our arbitrary definable set $Z$).

Instead, then, we have opted for the following notion. Suppose $Z$ is definable, and $\{X_t:t\in T\}$ is a definable family of subsets of $Z$ (all over $\emptyset$, say). We think of $\{X_t:t\in T\}$ as large `to order $k$' if for generic $t\in T$ and  $x_1,\dots,x_k\in X_t$ generic independent over $t$, the elements $x_1,\dots,x_k$ are also independent generics in $Z$ over $\emptyset$. One thinks of the $\{X_t\}$ as `filling' $Z$ to order $k$, so there are enough of the $X_t$ around to cover an independent $k$-tuple from $Z$. Correspondingly, it is easy to see that for fixed $Z$ and fixed $0<m<\dim Z$, the dimension of any family satisfying this condition (to order $k$) must tend to $\infty$ as $k$ does.

The main result of this subsection, then, answers the question posed above for our chosen notion of largeness. Namely, we show that for any definable set $Z$ and $0<m<\dim Z$, one can find a $k$-large (in the above sense) definable family of $m$-dimensional subsets of $Z$:

 \begin{proposition}\label{P: large families} Let $Z$ be a $\emptyset$-definable set of dimension $n\geq 2$. Then for all $k$ and $m$ with $0<m\leq n$ there exists a $\emptyset$-definable family $\{X_t: t\in T\}$ of subsets of $Z$ with the following properties: 
    \begin{enumerate}
        \item $T$ is stationary. 
        \item For generic $t\in T$ the set $X_t$ is a stationary $m$-dimensional subset of $Z$. 
        \item For generic $t\in T$ and $x_k,\dots, x_k\in X_t$ generic independent, $x_k,\dots, x_k$ are independent generics in $Z$ over $\emptyset$. 
    \end{enumerate}
\end{proposition}

\begin{proof} Our overarching strategy is as follows: first we show a related result in the case $S=M^2$ and $m=1$ (which implies the proposition in this case, but we don't need this); next, we use this initial result to prove an approximation to the proposition in the general case that $S=M^n$, with the only gap being that we sacrifice the stationarity of generic $X_t$. Finally, for arbitrary $S$ we set up a finite correspondence between $S$ and $M^n$, and then `transfer' our construction for $M^n$ through this correspondence. At this point we find it convenient to work in the language of stationary types; this will enable us in particular to fill the gap mentioned above; once we have the desired result for types, it is then easy to extract the desired definable family using compactness.

To begin, let us give our auxiliary (well known) result in the case $S=M^2$:
\begin{lemma}\label{good sets for M^2}
Let $\mathcal C=\{C_t:t\in T\}$ be a $\emptyset$-definable, faithful, generically strongly minimal family of plane curves in $M$, where $\dim(T)=k$. Let $t\in T$ be generic, and let $x_1,\dots,x_k\in C_t$ be independent generics over $t$. Then $t\in\acl(x_1,\dots,x_k)$.
\end{lemma}
\begin{proof} 
Note that $\dim(tx_1\dots x_k)=2k$, counting from left to right. Assume the lemma fails, and let $t'$ be an independent realization of $\tp(t/x_1\dots x_k)$ over $tx_1\dots x_k$. Then by assumption we get (a) $\dim(tt'x_1\dots x_k)>2k$, and (b) $\dim(t'/t)>0$. By (b) $t'\neq t$, so by the faithfulness of $\mathcal C$, $C_t\cap C_{t'}$ is finite, and thus $x_1,\dots,x_k\in\acl(tt')$. So by (a), $\dim(tt')>2k$, contradicting that $t,t'\in T$.

\end{proof}

We now use Lemma \ref{good sets for M^2} to show our approximate result for $M^n$:

\begin{lemma}\label{good sets for M^k}
    There is a $\emptyset$-definable family $\{Y_t: t\in T\}$ of subsets of $M^n$ with the the following properties: 
    \begin{enumerate}
        \item $T$ is stationary.
        \item $\dim(Y_t)=m$ for each $t\in T$.
        \item For any generic $t\in T$, and any independent generics $x_1,\dots, x_k\in Y_t$ over $t$, $x_1,\dots,x_k$ are also independent generics in $M^n$ over $\emptyset$. 
    \end{enumerate}
\end{lemma}
\begin{proof}
    Let $\{C_s:s\in S\}$ be a $\emptyset$-definable, faithful, generically strongly minimal family of plane curves, where $\dim(S)=k$. We may assume that for each $s$ the left projection $C_s\rightarrow M$ is finite-to-one. 
    
    Let $T=S^{n-m}$. Note that $T$ is stationary. Now we define the family $\{Y_t:t\in T\}$ as follows: for $t:=s_{m+1},\dots, s_n$, we let $$Y_t=\{(x_1,\dots,x_n):(x_m,x_j)\in C_{s_j}\textrm{ for }j=m+1,\dots,n\}.$$
    
    Using that each left projection $C_{S_j}\rightarrow M$ is finite-to-one, it is easy to see that the leftmost projection $Y_t\rightarrow M^m$ is finite-to-one with generic image in $M^m$; thus each $X_t$ has dimension $m$.
    
    Now let $t=(s_{m+1},\dots,s_n)\in T$ be generic, and let $x_1,\dots,x_k$ be independent generics of $Y_t$ over $t$. The key point is the following:
    
    \begin{claim} $t\in\acl(x_1,\dots,x_k)$.
    \end{claim}
    \begin{proof} For each $i\leq k$, let $x_i=(x_i^1,\dots,x_i^n)$, and set $y_i=(x_i^1,\dots,x_i^m)$. Note that each $x_i$ is interalgebraic with $y_i$ over $t$, by our remarks on the projection $Y_t\rightarrow M^m$.
    
    Now we have $\dim(x_1\dots x_k/t)=mk$, and thus also $\dim(y_1\dots y_k/t)=mk$. Since $y_1,\dots,y_k$ contains $mk$ total coordinates from $M$, these coordinates must all be independent generics over $t$. In particular, for each $j>m$, $x_1^m,\dots,x_k^m$ are independent generics in $M$ over $s_j$. This implies that $(x_1^m,x_1^j),\dots,(x_k^m,x_k^j)$ are independent generics in $C_{s_j}$ over $s_j$. Clearly $s_j$ is generic in $S$, because $t$ is generic in $T$; thus by Lemma \ref{good sets for M^2}, we get $s_j\in\acl(x_1^m,x_1^j,\dots,x_k^m,x_k^j)$. Combining for all $j$ then gives the claim.
    \end{proof}
    
    Now by the claim we have $$\dim(x_1\dots x_k)=\dim(tx_1\dots x_k)$$ $$=\dim(t)+\dim(x_1\dots x_k/t)=k(n-m)+km=kn.$$ So $x_1,\dots,x_k$ are independent generics in $M^n$, as desired.
\end{proof}

\begin{remark} Proposition \ref{good sets for M^k} is obviously false if we allow the case $m=0$. There is only one subtle place in the proof of Lemma \ref{good sets for M^k} where $m>0$ is used: indeed, the reader may wish to verify that the definition of the $Y_t$ from the $C_s$ does not make any sense if $m=0$.
\end{remark}

We now return to the general case of the proposition. By weak elimination of imaginaries in $\mathcal M$, every generic of $Z$ is interalgebraic with a real tuple. In particular, extracting a basis from such a tuple, we get that every generic of $Z$ is interalgebraic with a generic of $M^n$. Since $M^n$ is stationary, the reverse also holds, i.e. every generic of $M^n$ is interalgebraic with a generic of $Z$.

Now let $\{Y_t:t\in T\}$ be as in the lemma. Note that we may assume $k\geq 1$. Let $t_0\in T$ be generic, and let $y\in Y_{t_0}$ be generic over $t_0$. Then $y$ is generic in $M^n$ over $\emptyset$. Let $x$ be a generic of $Z$ which is interalgebraic with $y$, and let $p=\operatorname{stp}(x/t_0)$. So $\dim(p)=\dim(x/t_0)=\dim(y/t_0)=m$.

\begin{lemma} If $x_1,\dots,x_k$ are independent realizations of $p$ over $t_0$, then $x_1,\dots,x_k$ are independent generics of $Z$ over $\emptyset$.
\end{lemma}
\begin{proof} Since each $\tp(x_i/t_0)=\tp(x/t_0)$, for each $i$ we can find a generic $y_i$ of $Y_{t_0}$ which is interalgebraic with $x_i$. Since the $x_i$ are independent, we have $\dim(x_1,\dots,x_k/t_0)=mk$, and thus by interalgebraicity $\dim(y_1,\dots,y_k/t_0)=mk$. Thus $y_1,\dots,y_k$ are independent generics of $Y_{t_0}$ over $t$, so by assumption $y_1,\dots,y_k$ are also independent generics of $M^n$ over $\emptyset$. Then $\dim(y_1,\dots,y_k)=nk$, so by interalgebraicity again $\dim(x_1,\dots,x_k)=nk$, and thus $x_1,\dots,x_k$ are independent generics of $Z$.
\end{proof}

Now by compactness, there is a $\emptyset$-definable family $\{X_t:t\in T'\}$ of $m$-dimensional subsets of $Z$, where $T'$ is a generic subset of $T$, and such that $X_{t_0}$ is stationary with generic type $p$. Then this family clearly satisfies (1)-(3) in the proposition statement, so we are done.

\end{proof}

\subsection{Sets Transverse to Any Family} 

The result of this subsection was inspired by the following question: given a definable family $\{X_s:s\in S\}$ of codimension $r$ subsets of a definable set $Z$, can one always find an $r$-dimensional definable set $Y\subset Z$ having finite intersection with each $X_s$? This type of question arises, for example, when trying to extract a 2-dimensional very ample family of plane curves from a higher dimensional such family (see Proposition \ref{P: down}).  

In order to answer the question above, it is convenient to study a more general question. Namely, for definable sets $X,Y\subset Z$, say that $X$ and $Y$ are \textit{transverse in} $Z$ if $\dim(X\cap Y)\leq\dim X+\dim Y-\dim Z$. We wish to investigate whether, given a definable family $\{X_s:s\in S\}$ of subsets of a definable set $Z$, one can find a definable $Y\subset Z$ of any prescribed dimension which is transverse to each $X_s$ in $Z$; and more generally whether one can find a large family of such sets (in the sense of Proposition \ref{P: large families}). 

As stated, it is not always possible to find such $Y$: indeed, if the prescribed dimension of $Y$ is small enough (but still at least $0$), transversality could require $Y$ to be disjoint from each $X_s$, which is clearly impossible if the $X_s$ cover $Z$. So at the least, we must allow $Y$ to be non-empty. Conversely, Proposition \ref{P: transverse families} shows that this is the only obstacle in general. Namely, we show that given $\{X_s:s\in S\}$ inside $Z$, one can always find a large family $\{Y_t:t\in T\}$ inside $Z$ so that each $X_s$ and $Y_t$ are either transverse in $Z$, or $\dim(X_s)+\dim(Y_t)<\dim Z$ and $X_s\cap Y_t$ is finite:

\begin{proposition}\label{P: transverse families} Fix $0\leq i\leq m$ and $0<j<m$. Let $Z$ be a definable set of dimension $m$, and let $\{X_s:s\in S\}$ be a definable family of subsets of $Z$ with each $X_s$ of dimension at most $i$. Then for all $k$ there are a non-empty definable set $T$ and a definable family $\{Y_t:t\in T\}$ of $j$-dimensional subsets of $Z$ satisfying the following:
\begin{enumerate}
    \item For each $s\in S$ and each $t\in T$ we have $$\dim(X_s\cap Y_t)\leq\max\{0,i+j-m\}$$ (here we take $\dim(\emptyset)=-\infty$ by convention).
    \item For any generic $t\in T$ and any independent generic $z_1,\dots,z_k\in Y_t$, $z_1,\dots,z_k$ are independent generics in $Z$.
\end{enumerate}
\end{proposition}
\begin{proof} Throughout the proof, we will assume $k\geq 2$. This is harmless, as higher $k$-values are even more restrictive. We also assume $\{X_s:s\in S\}$ and $Z$ are $\emptyset$-definable.

We will induct on $j$. First suppose $j=1$. If $i=m$ then (1) holds automatically; so we only need to ensure (2), and this is possible by Proposition \ref{P: large families}.

So suppose $i\leq m-1$.

\begin{claim} There is some $n$ such that for any independent generics $z_1,\dots,z_n\in Z$, there does not exist $s\in S$ with $z_1,\dots,z_n\in X_s$.
\end{claim}
\begin{proof} Suppose $z_1,\dots,z_n\in Z$ are independent generics in $Z$, and there is some $s\in S$ with $z_1,\dots,z_n\in X_s$. Then $\dim(z_l/s)\leq i\leq m-1$ for each $l$, and so $\dim(s,z_1\dots z_n)\leq \dim(S)+n(m-1)$. But by assumption $\dim(z_1\dots z_n)=mn$. So $mn\leq m+n(m-1)$, which simplifies to $n\leq \dim(S)$. That is, choosing any $n>\dim(S)$  the claim  holds.
\end{proof}

Now let $n$ be as in the claim. Without loss of generality, we assume $n=k$, since we are allowed to increase either quantity until it matches the other. By Proposition \ref{P: large families} there is a definable family $\{Y_t:t\in T\}$ of 1-dimensional subsets of $Z$ satisfying (2). Without loss of generality, assume this family is $\emptyset$-definable.

\begin{claim} For $t\in T$ generic and $s\in S$ arbitrary, $X_s\cap Y_t$ is finite (possibly empty).
\end{claim}
\begin{proof}
Suppose not. Since $\dim(Y_t)=1$, we can  find independent generics $z_1,\dots,z_n$ of $Y_t$, each of which belongs to $X_s$. By (2), $z_1,\dots,z_n$ are independent generics in $Z$ and all belong to $X_s$, which contradicts the choice of $n$.
\end{proof}

Finally, by compactness (1) holds after removing a non-generic definable subset of $T$, and this process clearly does not affect (2). We have thus completed the case $j=1$.

Now assume $j>1$ and the proposition holds for all $j'<j$. Let $\{D_u:u\in U\}$ satisfy (1) and (2) with $j$ replaced by $j-1$ (and all other data unchanged).

\begin{lemma}\label{L: first} We may assume that for each $z\in Z$, the set $\{u\in U:z\in D_u\}$ has dimension at most $\dim U+j-1-m$.
\end{lemma}
\begin{proof}
Without loss of generality, $\{D_u:u\in U\}$ is $\emptyset$-definable. Note that $\{D_u:u\in U\}$ still satisfies (1) and (2) if we delete a non-generic portion of $U$ and a non-generic portion of each remaining $D_u$. So by compactness, it suffices to show that for generic $u\in U$ and generic $z\in D_u$ we have $\dim(u/z)=\dim U+j-1-m$. Now by genericity we have $\dim(u)=\dim U$ and $\dim(z/u)=j-1$; while by (2) we have $\dim(z)=m$. The desired equality then follows by additivity.
\end{proof}
Continuing with the proof of the proposition, for $s\in S$ let us denote $$E_s:=\{u\in U:X_s\cap D_u\neq\emptyset\}.$$ We now view $\{E_s:s\in S\}$ as a definable family of subsets of $U$. Let $m'=\dim(U)$, and let $i'$ be the largest dimension of any $E_s$. Then by the $j=1$ case, we can find a family $\{W_t:t\in T\}$ of 1-dimensional subsets of $U$ satisfying (1) and (2) while replacing $Z$ with $U$, $m$ with $m'$, $i$ with $i'$, $j$ with 1, $\{X_t\}$ with $\{E_s\}$, and keeping the same $k$.

Finally, for $t\in T$ let $Y_t$ be the union of all $D_u$ for $u\in W_t$. Clearly $\{Y_t\}$ is a definable family of subsets of $Z$. As in the $j=1$ case, we will show that $\{Y_t\}$ satisfies (1) and (2) generically, and deduce that (1) and (2) can be satisfied fully after removing a small subset of $T$.

Absorbing all necessary parameters into the language, we assume henceforth, assume that all data discussed to this point is $\emptyset$-definable.

\begin{lemma} Let $t_0\in T$ be generic. Then $\dim(Y_{t_0})=j$.
\end{lemma}
\begin{proof}
By definition $Y_{t_0}$ is the union of a 1-dimensional family of $j-1$-dimensional sets, so it is clear that $\dim(Y_{t_0})\leq j$. Now let $u_0\in W_{t_0}$ be generic over $t_0$, and let $z_0\in D_{u_0}$ be generic over $t_0u_0$. So $z_0\in Y_{t_0}$. Applying (2) to the family $\{W_t\}$, we get that $u_0$ is generic in $U$; then applying (2) again, this time to $\{D_u\}$, we get that $z_0$ is generic in $Z$. Now $\dim(u_0/t_0)=1$ and $\dim(z_0/t_0u_0)=j-1$, so $\dim(u_0z_0/t_0)=j$. If $u_0\in\acl(t_0z_0)$ then $\dim(z_0/t_0)=j$, which implies the lemma. So assume toward a contradiction that $u_0\notin\acl(t_0z_0)$. Since $u_0\in W_{t_0}$ and $\dim(W_{t_0})=1$, $u_0$ is generic in $W_{t_0}$ over $t_0z_0$. By (2) again (applied to $\{W_t\}$ after absorbing the parameter $z_0$), this implies that $u_0$ is generic in $U$ over $z_0$. Thus $u_0$ and $z_0$ are independent, which gives that $z_0$ is generic in $Z$ over $u_0$. But this is not true, because $z_0\in D_{u_0}$ and $\dim(D_{u_0})=j-1<m$.
\end{proof}

We now proceed to verify a generic version of (1): 

\begin{lemma}\label{L: transverse} Let $t_0\in T$ be generic, and let $s_0\in S$ be arbitrary. Then $$\dim(X_{s_0}\cap Y_{t_0})\leq\max\{0,i+j-m\}.$$
\end{lemma}
\begin{proof}
First assume $i+j>m$. Then $i+(j-1)-m\geq 0$, so by (1) applied to the family $\{D_u\}$ the intersection of $X_{s_0}$ with any $D_u$ has dimension at most $i+(j-1)-m$. Since $Y_{t_0}$ is the union of a 1-dimensional family of the $D_u$,  necessarily $X_{s_0}$ and $Y_{t_0}$ intersect in dimension at most $(i+(j-1)-m)+1=i+j-m$, and we are done.

Now assume $i+j\leq m$. So by (1) (applied to $\{D_u\}$), the intersection of $X_{s_0}$ with each $D_u$ is finite. Recalling the definition of $E_s$ above, we first show:
\begin{claim}
For each $s\in S$, $\dim(E_s)\leq\dim(U)-1$.
\end{claim}
\begin{proof}
Recall that we defined $m'=\dim (U)$. Let $u\in E_s$ be arbitrary. It will suffice to show that $\dim(u/s)\leq m'-1$. Now since $u\in E_s$, there is some $z\in X_s\cap D_u$. Since $z\in X_s$, $\dim(z/s)\leq i$. Since $z\in D_u$ and by Lemma \ref{L: first}, $\dim(u/z)\leq m'+j-1-m$. So by additivity, $\dim(uz/s)\leq m'+j-1-m+i$. Since $i+j\leq m$, we simplify to $\dim(uz/s)\leq m'-1$. In particular, $\dim(u/s)\leq m'-1$, as desired.
\end{proof}
Now by the claim and the definition of $i'$ and $m'$, we get $i'\leq m'-1$. So in calling on the $j=1$ case to obtain the family $\{W_t\}$, we had $i'+(1)-m'\leq 0$. Thus, by (2) (applied to the family $\{W_t\}$ again), the intersection of $E_{s_0}$ and $W_{t_0}$ is also finite. Now by definition, $X_{s_0}\cap Y_{t_0}$ is the union of the $X_{s_0}\cap D_u$ for $u\in W_{t_0}$; and by definition of $E_{s_0}$, it suffices to consider only those $u\in E_{s_0}$. Then $X_{s_0}\cap Y_{t_0}$ is the union of finitely many finite sets, so has dimension at most 0. This completes the proof of the lemma. 
\end{proof}

Now we verify a generic version of (2): 

\begin{lemma}\label{L: transverse family is large} Let $t\in T$ be generic, and let $z_1,\dots,z_k\in Y_t$ be independent generics over $t$. Then $z_1,\dots,z_k$ are independent generics in $Z$.
\end{lemma}
\begin{proof}
By definition, there are $u_1,\dots,u_k\in W_t$ with each $z_l\in D_{u_l}$ for $l=1,\dots,k$. Now by assumption we have $\dim(z_1\dots z_k/t)=jk$. Since each $z_l\in D_{u_l}$, also  $\dim(z_1\dots z_k/tu_1\dots u_k)\leq(j-1)k$. So $\dim(u_1\dots u_k/t)\geq k$. In particular, this implies that $u_1,\dots,u_k$ are independent generics in $W_t$ over $t$. By (2), $u_1,\dots,u_k$ are also independent generics in $U$ over $\emptyset$. We conclude the following two claims:
\begin{claim} $\dim(t/u_1\dots u_k)=\dim(T)-k(\dim (U)-1)$.
\end{claim}
\begin{proof} By what we have seen above,  we have $\dim(tu_1\dots u_k)=\dim T+k$, and $\dim(u_1\dots u_k)=k\cdot\dim(U)$. Thus $\dim(t/u_1\dots u_k)=\dim(T)+k-k\cdot\dim(U)$, which is the same as the desired quantity.
\end{proof}
\begin{claim} $\dim(u_1\dots u_k/z_1\dots z_k)\leq k(\dim U+j-1-m)$.
\end{claim}
\begin{proof} By Lemma \ref{L: first} each $u_l$ has dimension at most $\dim(U)+j-1-m$ over $z_l$. The claim is now obvious.
\end{proof}
By the two claims and additivity, $\dim(tu_1\dots u_k/z_1\dots z_k)\leq\dim T+k(j-m)$; while by assumption $\dim(tu_1\dots u_kz_1\dots z_k)\geq\dim(tz_1\dots z_k)=\dim T+jk$. So by additivity, $\dim(z_1\dots z_k)\geq km$, which shows that $z_1,\dots,z_k$ are independent generics in $Z$.
\end{proof}
Finally, by Lemma \ref{L: transverse} and compactness, after deleting a non-generic set $V\subset T$ the equality in (1) holds for all $t\in T-V$, and the conclusion of Lemma \ref{L: transverse family is large} is clearly unchanged. This completes the inductive step.
\end{proof}

\section{Very Ample Families of a Prescribed Dimension}
As we have seen, very ampleness in strongly minimal structures is equivalent to the existence of very ample families of plane curves. In applications, it may be convenient to have access to such families of a prescribed dimension (as is the case with non-local modularity, where the existence of any $k\geq 2$-dimensional family of plane curves implies the existence of families of all dimensions). In the present section, we show that this can always be achieved. Namely, we first use Proposition \ref{P: transverse families} to show that, given a $k$-dimensional very ample family of plane curves, we can extract an $l$-dimensional very ample sub-family for any $2\le l\le k$. We then adapt the usual proof from the non-locally modular setting (using composition and a field configuration) to construct very ample families of arbitrarily large dimension. While the general proof structure is the same as in the non-locally modular case, the need to preserve very ampleness will make each stage of the argument more delicate. In particular, once a field (say $K$) is interpreted in our set (say $X$), some care is needed  to examine the precise relationship between $X$ and $K$.

Once we have achieved the result described above, we end the section with a short application, in which (as mentioned in the introduction to the paper) we characterize very ampleness in terms of definable pseudoplanes. 

\subsection{Going Down}

Here we prove that very ample families of plane curves admit very ample subfamilies of any prescribed dimension. The proof amounts to an application of Proposition \ref{P: transverse families}:

\begin{proposition}\label{P: down} Let $X$ be strongly minimal and $A$-definable. Let $C$ be a strongly minimal plane curve in $X$, let $c=\operatorname{Cb}(C)$, and let $k=\dim(c/A)$. Assume that $C$ is very ample over $A$. Then for all $2\leq j<k$, there is a set $B\supset A$ such that $\dim(c/B)=j$ and $C$ is very ample over $B$.
\end{proposition}
\begin{proof}
Without loss of generality, assume $A$ is algebraically closed, as this does not affect any of the relevant computations. By Corollary \ref{C: curves to families} there is an $A$-definable very ample family $\{C_z:z\in Z\}$ of plane curves in $X$, and a generic $z_0\in Z$ over $A$, such that $C$ and $C_{z_0}$ have finite symmetric difference. Since $A$ is algebraically closed, we may assume $Z$ is stationary. If we apply an appropriate quotient to $Z$ (and potentially remove a small portion of $Z$) we may moreover assume $\{C_z\}$ is faithful, and thus $c$ and $z_0$ are interdefinable over $A$.

Let $m=\dim(Z)$, and let $S=\{(x,y)\in M^2\times M^2:x\neq y\}$, the set of pairs of distinct points in the plane. For $s=(x,y)\in S$ let $X_s=\{z\in Z:x,y\in C_z\}$; then very ampleness gives $\dim(X_s)\leq m-2$ for each $s$. Now let $i=m-2$ and $k=1$, and apply Proposition \ref{P: transverse families} to extract a family $\{Y_t:t\in T\}$ of $j$-dimensional subsets of $Z$. Let $A'\supset A$ so that $\{Y_t\}$ is $A'$-definable. Replacing $\{Y_t\}$ by its image under an appropriate automorphism, we may assume that $z_0$ is independent of $A'$ over $A$, so still generic in $Z$.

Let $t$ be generic in $T$ over $A'$, let $B=A't$, and let $z$ be generic in $Y_t$ over $B$. So $\dim(z/B)=\dim(Y_t)=j$. 

\begin{claim} $C_z$ is very ample over $B$.
\end{claim}
\begin{proof} Let $x\neq y$ be generics in $C_z$ over $B$. By Lemma \ref{L: Ind on parms} and the above, it suffices to show that $\dim(z/Bxy)\leq j-2$. Now since $x\neq y$, the pair $s=(x,y)$ belongs to $S$. Thus $z\in X_s\cap Y_t$. By (1) of Proposition \ref{P: transverse families} $\dim(X_s\cap Y_t)\leq i+j-m$, and since $i=m-2$ this simplifies to $j-2$. Since $X_s\cap Y_t$ is $Bxy$-definable, the desired conclusion follows.
\end{proof}

Finally, by (2) of Proposition \ref{P: transverse families}, note that $z$ is generic in $Z$ over $A'$. So since $Z$ is stationary, there is an automorphism over $A'$ sending $z$ to $z_0$. Let $B_0$ be the image of $B$ under this automorphism; then clearly $C_{z_0}$ is very ample over $B_0$, and $\dim(z_0/B_0)=\dim(z/B)=j$. Since $c$ and $z_0$ are interdefinable over $B$, $\dim(c/B_0)$ is also $j$, which completes the proof.

\end{proof}

\subsection{Going Up} 
Our next goal is to show that if $X$ is a very ample strongly minimal set then $X$ has arbitrarily large very ample families of plane curves. Before proceeding, we need the following lemma, that is by now standard:
\begin{lemma}\label{L: composing up}
Let $X$ be a $A$-definable strongly minimal set, and let $C$, $D$, and $E\subset D\circ C$ be strongly minimal plane curves in $X$, with canonical bases $c$, $d$, and $e$ respectively. Assume that $c\ind_Ad$, and $\dim(e/A)\leq\dim(c/A)=\dim(d/A)=k$ for some $k\geq 2$. Then $k\in \{2,3\}$, and there is a definable strongly minimal expansion of an algebraically closed field which is internal to $X$. 
\end{lemma}
\begin{proof}
 Since the proof is standard, we will be brief. We may assume $A=\emptyset$. Clearly $\dim(e)\ge k$ (see also\cite[Lemma 3.18]{ElHaPe}). So we only have to deal with the case that $\dim(e)=k$. In this case, take $(x,z)\in E$ generic and $y$ such that $(x,y)\in C$ and $(y,z)\in D$. Then $\{c,d,e, x,y,z\}$ is a group configuration (see \cite[Lemma 3.20]{ElHaPe} for the details). By a result of Hrushovski (see \cite[Theorem 3.27]{PoiGroups}) this implies that $k\le 3$, and if $\dim(k)\geq 2$ (which holds by assumption), a strongly minimal field is definable in $X^{\textrm{eq}}$. 
\end{proof}
We now give the main result of this subsection:
\begin{proposition} \label{P: going up}
Suppose $X$ is strongly minimal, $A$-definable, and very ample. Then for any $k\geq 2$ one can find a set $B\supset A$, and a strongly minimal plane curve $C$ in $X$, such that $\dim(\operatorname{Cb}(C)/B)\geq k$ and $C$ is very ample in $X$ over $B$.
\end{proposition}
\begin{proof}
    For simplicity of notation, we assume that $A=\0$. Assume the proposition fails. By very ampleness, there is a strongly minimal plane curve $C$ in $X$ (with canonical base $c$, say) which is very ample in $X$ over some set $B$. Fix such $(C,c,B)$ so that the value $k:=\dim(c/B)$ is maximal (that this is possible follows by the assumed failure of the proposition). We may assume moving forward that $B=\0$. Note also that $k\geq 2$ by Lemma \ref{L: very ample implies nlm}.

    Let $d\models\operatorname{tp}(c/B)$ with $c\ind d$. So $\dim(d)=k$, and $d=\operatorname{Cb}(D)$ for some strongly minimal plane curve $D$ in $X$ which is also very ample in $X$ over $\0$. Let $E\subset D\circ C$ be strongly minimal. Then, by Lemma \ref{L: composition} and Corollary \ref{C: VA components}, $E$ is also very ample over $\0$. So by the choice of $k$, letting $e=\operatorname{Cb}(E)$, we have $\dim(e)\leq k$. Thus, by Lemma \ref{L: composing up}, $k\in\{2,3\}$ and there is a strongly minimal expansion $(K,+,\cdot,\dots)$ of an algebraically closed field internal to $X$.

    Since $K$ is internal to $X$, clearly $X$ and $K$ are non-orthogonal. So by Corollary \ref{C: field va int equivalence}, there is a definable embedding of $X$ into some $K^n$. Let us choose, among all definable embeddings of cofinite subsets of $X$ into any $K^n$, one which minimizes the value $n$. So, replacing $X$ with a cofinite subset if necessary, we assume that $X\subset K^n$ and no cofinite subset of $X$ embeds into any smaller power of $K$. Without loss of generality, we will assume that $K$ and $X$ are $\emptyset$-definable.
   
   In what follows, it will be convenient to introduce the following notation: if $S\sub K^m$ is definable, we let $S^{\mathrm{aff}}$ be the intersection of all affine linear subspaces of $K^n$ which contain a large subset of $S$. It is easy to see that $S^{\mathrm{aff}}$ is itself affine linear and contains a large subset of $S$, and moreover that $S^{\mathrm{aff}}$ is definable over any set of parameters defining $S$. Note also that by the choice of $n$ we have $X^{\mathrm{aff}}=K^n$ (since any proper affine linear subspace of $K^n$ is in definable bijection with a lower power of $K$).
   
   Let $\CH=\{H_t:t\in T\}$ be the family of hyperplanes in $K^{2n}$. So $\dim(T)=2n$. We note the following basic observation, for future reference: 
   \begin{claim}\label{C: dim of hyperplanes through tuple}
        If $D\sub K^{2n}$ is a $d$-dimensional affine linear subspace, then $\dim(\{t\in T: D\sub H_t\})=2n-d-1$. 
   \end{claim}
   
   Now let $t\in T$ be generic. Then by Proposition \ref{P: field va}, $H_t\cap X^2$ is a very ample (over $\emptyset$) plane curve in $X$. Let $S\subset H_t\cap X^2$ be strongly minimal. Then by Corollary \ref{C: VA components}, $S$ is also very ample in $X$ over $\emptyset$. So by the choice of $k$, letting $s=\operatorname{Cb}(S)$ we have $\dim(s)\leq k\leq 3$, and thus $\dim(s)\leq 3$.

   Since $S$ is a strongly minimal component of $H_t\cap X^2$, clearly $s\in\acl(t)$. So $\dim(ts)=\dim(t)=2n$, and thus $\dim(t/s)\geq 2n-3$. On the other hand, by definition $H_t$ contains $S$, so since $H_t$ is affine linear it must also contain $S^{\mathrm{aff}}$. Let $d:=\dim(S^{\mathrm{aff}})$; then by Claim \ref{C: dim of hyperplanes through tuple} we have $\dim(t/s)\leq 2n-d-1$. So, combining with $\dim(t/s)\geq 2n-3$, we get $2n-3\leq 2n-d-1$, and thus $d\leq 2$. Note also that if equality holds (i.e. $d=2$), we have $2n-3\leq\dim(t/s)\leq 2n-3$, so that $\dim(t/s)=2n-3$, and thus $\dim(s)=3$.\\

   Now the main point of the proof is the following:

   \begin{lemma}\label{L: definable field} There is a definable strongly minimal algebraically closed field $(F,+_F,\cdot_F)$, whose underlying set $F$ almost coincides with $X$.
   \end{lemma}
   \begin{proof} First note that if $\dim(X^{\mathrm{aff}})=1$, then $X$ almost coincides with $X^{\mathrm{aff}}$, and $X^{\mathrm{aff}}$ is in definable bijection with $K$; in particular, this is enough to prove the lemma (setting $F=X^{\mathrm{aff}}$ with the field structure inherited from $K$). So we may assume that $n=\dim(X^{\mathrm{aff}})\geq 2$, and thus (as above) $\dim(s)=3$. Our goal is to construct a definable transitive group action of a 3-dimensional group on a set almost equal to $X$, and apply Hrushovski's classification (\cite[Theorem 3.27]{PoiGroups}) of such actions. We first proceed with the following two claims, which give us the group we will use:
   
   \begin{claim}\label{Cl: S covers K^2} Let $\pi:K^{2n}\rightarrow K^n$ be either the leftmost or rightmost projection. Then the restriction of $\pi$ to $S^{\mathrm{aff}}$ is surjective.
   \end{claim}
   \begin{proof} By Remark \ref{R: non trivial}, $\pi$ is finite-to-one on $S$, and thus $\pi(S)$ contains a cofinite subset of $X$. By definition, $S^{\mathrm{aff}}$ contains a cofinite subset of $S$, so clearly $\pi(S^{\mathrm{aff}})$ contains a cofinite subset of $\pi(S)$. In particular, $\pi(S^{\mathrm{aff}})$ also contains a cofinite subset of $X$, and thus contains $X^{\mathrm{aff}}=K^n$ (because projections of affine sets are affine).
   \end{proof}
   
   We next use Claim \ref{Cl: S covers K^2} to show the following:
   
   \begin{claim}\label{Cl: S^aff linear} $S^{\mathrm{aff}}$ is the graph of an invertible affine linear map $L:K^2\rightarrow K^2$.
   \end{claim} 
   \begin{proof} By Claim \ref{Cl: S covers K^2} we have $d\geq n\geq 2$; so since we already had $d\leq 2$, we get $d=n=2$. Then by Claim \ref{Cl: S covers K^2} again, each of the two projections $S^{\mathrm{aff}}\rightarrow K^2$ is a surjective linear map between two-dimensional affine linear spaces; thus each of these projections is a linear bijection, which implies the claim.
   \end{proof}
   
   Let $L:K^2\rightarrow K^2$ be as in Claim \ref{Cl: S^aff linear}. Since $S$ projects almost onto $X$ in both directions, and by the strong minimality of $X$, it follows that $L(X)$ almost coincides with $X$. By the strong minimality of $S$, it moreover follows that $S$ almost coincides with the graph of the restriction of $L$ to $X$. In particular, it follows easily that $s$ is interdefinable with the canonical parameter of $L$.

   Now let $G$ be the generic stabilizer of $X$ in $\mathrm{AGL}_2(K)$ -- that is, the set of all invertible affine linear maps $g:K^2\rightarrow K^2$ such that $g(X)$ almost coincides with $X$. Clearly, $G$ is a $\emptyset$-definable subgroup of $\mathrm{AGL}_2(K)$. By the above remarks (since $s$ is definable from $L$), we have $\dim(G)\geq\dim(s)=3$. Let $G^0$ be the connected component of $G$; so $\dim(G^0)\geq 3$ and $G^0$ is $\emptyset$-definable.

   We now have the group for our action; our next goal is to find an appropriate set which is acted on transitively. Now $G^0$ comes equipped with both an action on $K^2$, and a generic action on $X$: so if $g\in G^0$ and $x\in X$ then $g(x)\in K^2$ is always defined, while if $g$ and $x$ are independent generics then $g(x)$ is also a generic of $X$. Fix $x_0\in X$ generic and independent from $s$, and let $G^0(x_0)$ be the orbit of $x_0$ under all of $G^0$. Note, then, that $G^0(x_0)$ is definable over $x_0$. Our aim is to show that $G^0(x_0)$ almost coincides with $X$, and subsequently build a field structure on $G^0(x_0)$. We will do this via a sequence of claims. To start, we show:

   \begin{claim}\label{Cl: action is transitive} $G^0(x_0)$ almost contains $X$.
   \end{claim}
   \begin{proof}
   Since $x_0$ is generic over $s$, there is $y_0\in X$ such that $(x_0,y_0)$ is generic in $S$. So $y_0=L(x_0)$. Now by genericity we have $\dim(x_0y_0/s)=1$; Note then that $\dim(x_0y_0/\emptyset)$ cannot also be 1, because then we would have $s\in\acl(\emptyset)$, contradicting that $\dim(s)=3$. So $\dim(x_0y_0)=2$. In particular, $y_0$ is generic in $X$ over $x_0$. Since $y_0=L(x_0)$ we see that $y_0\in G^0(x_0)$, and the claim follows.
   \end{proof}

   Using Claim \ref{Cl: action is transitive}, we now show:

   \begin{claim} $G^0(x_0)$ is one-dimensional, and $X$ is (up to a finite set) one of its strongly minimal components.
   \end{claim}
   \begin{proof} Let $D$ be the set of $g\in G^0$ with $g(x_0)\in X$. Then $D$ is definable over $x_0$ and contains all generics of $G^0$ over $x_0$, which implies that $D$ is generic in $G^0$. So there are $g_1,\dots,g_m\in G^0$ such that the translates $g_iD$ cover $G^0$. Then $G^0(x_0)$ is contained in the union of the $g_iD(x_0)$, which by definition is contained in the union of the $g_i(X)$. Now each $g_i(X)$ is strongly minimal, which implies that $G^0(x_0)$ is contained in a one-dimensional definable set. The claim then follows by Claim \ref{Cl: action is transitive}.
   \end{proof}

   Finally, we show:

   \begin{claim}\label{Cl: sm orbit} $G^0(x_0)$ is strongly minimal, and thus almost coincides with $X$.
   \end{claim}
   \begin{proof} As a transitive $G^0$-set, $G^0(x_0)$ is isomorphic to a coset space $G/H$, where $H$ is the stabilizer of $x_0$. It is easy to see that this identification is definable (that is, the usual identification from basic algebra is clearly definable). In this case, the strong minimality of $G^0/H$ follows easily from the fact that $G^0$ is connected: indeed, given a partition of $G^0/H$ into disjoint definable infinite sets $Z_1$ and $Z_2$, the preimages of the $Z_i$ in $G^0$ would partition $G^0$ into two generic subsets, contradicting connectedness.
   \end{proof}

   Now by Claim \ref{Cl: sm orbit}, there is a definable (transitive by definition) action of $G^0$ on the strongly minimal set $G^0(x_0)$, and $\dim(G^0)\geq 3$. So by the classification of such actions (\cite[Theorem 3.27]{PoiGroups}), $\dim(G^0)=3$ and $G^0(x_0)$ is in definable bijection with $\mathbb P^1(F)$ for some definable, strongly minimal, algebraically closed field $F$. Deleting a point, we then recover $\mathbb A^1(F)=F$; then by Claim \ref{Cl: sm orbit} again, $X$ is (up to a definable bijection) almost equal to $F$, proving Lemma \ref{L: definable field}.
   \end{proof}
   
Now by Lemma \ref{L: definable field}, we may assume after changing finitely many points that there is a definable field structure on the set $X$. It is then an easy exercise to find arbitrarily large very ample families of plane curves (for example, the family of degree $d$ polynomial maps $X\rightarrow X$ is very ample for each $d$). In particular, there are very ample (over $\emptyset$) strongly minimal plane curves $C\subset X^2$ with $\dim(\operatorname{Cb}(C))>k$, contradicting our initial assumption and thus proving Proposition \ref{P: going up}.   
\end{proof}

Finally, we now deduce:

\begin{corollary}\label{C: families of all dims}
    Let $X$ be a very ample strongly minimal set. Then for every $k\geq 2$ there is a faithful, very ample family $\{C_t:t\in T\}$ of plane curves in $X$ where $\dim(T)=k$.
\end{corollary}
\begin{proof} Assume $X$ is $A$-definable. Applying Proposition \ref{P: going up} and then Proposition \ref{P: down}, we can find a set $B\supset A$ and a strongly minimal plane curve $C\subset X^2$ which is very ample in $X$ over $A$ and satisfies $\dim(\operatorname{Cb}(C)/B)=k$. After editing finitely many points of $C$ and applying compactness, there is then a $B$-definable faithful family of plane curves in $X$ whose generic members are the conjugates of $C$ over $B$. The result then follows by applying Proposition \ref{P: family VA iff curve VA} to this family.
\end{proof}

\subsection{Very Ampleness and Pseudoplanes} We now show, as described in the introduction, that a strongly minimal set $X$ is very ample (i.e., admits a very ample plane curve) if and only if there is a definable pseudoplane on (a large subset of) $X^2$. Let us begin by formally recalling the definition:

\begin{definition}\label{D: pseudoplane}
    A \textit{pseudoplane} consists of sets $P$ and $L$, and an incidence relation $I\subset P\times L$, satisfying the following:
    \begin{enumerate}
         \item For each $l\in L$ there are infinitely many $p\in P$ with $(p,l)\in I$.
        \item For each $p\in P$ there are infinitely many $l\in L$ with $(p,l)\in I$.
        \item For each $l\neq l'\in L$ there are only finitely many $p\in P$ with $(p,l),(p,l')\in I$.
        \item For each $p\neq p'\in P$ there are only finitely many $l\in L$ with $(p,l),(p',l)\in I$.
    \end{enumerate}
\end{definition}

One thinks of $P$ as the set of `points' of a plane, and $L$ as the set of `lines.' As mentioned above, already Zilber in his thesis proved the "weak trichotomy theorem", asserting that  a strongly minimal structure is not locally modular if and only if it defines a pseudoplane. However, if one wishes to make the natural identification $P:=M^2$ (or a large subset thereof), one can only guarantee conditions (1)-(3) in the definition. We now point out that the missing data in order to identify $P$ (generically) with $M^2$ is precisely very ampleness.

\begin{proposition}\label{P: VA = nice PP} Let $X$ be strongly minimal. Then the following are equivalent:
\begin{enumerate} 
\item $X$ is very ample.
\item There is a definable pseudoplane $(P,L,I)$ where $P$ is a generic subset of $X^2$.
\end{enumerate}
\end{proposition}
\begin{proof} First, assume $X$ is very ample. Then by Proposition \ref{P: down}, there is a definable faithful very ample family of plane curves $\mathcal C=\{C_t:t\in T\}$ in $X$, where $\dim(T)=2$. Let $C\subset X^2\times T$ be the graph of $\mathcal C$. Let $P$ be the set of $x\in X^2$ such that the set $C^x\subset T$ has dimension 1. Let $L$ be the set of $t\in T$ such that $C_t\cap P$ is infinite. Let $I=C\cap(P\times L)$. It is easy to check that $P$ is generic in $M^2$ and $L$ is cofinite in $T$. It then follows easily that axioms (1) and (2) of Definition \ref{D: pseudoplane} hold for $(P,L,I)$ (these computations are essentially carried out in section 2 of \cite{CasACF0}). To conclude, we note that (3) is a restatement of the faithfulness of $\mathcal C$, and (4) is a restatement of the very ampleness of $\mathcal C$.

Now assume $(P,L,I)$ is a definable pseudoplane, where $P$ is generic in $X^2$. We will show that $\{I_l:l\in L\}$ is a two-dimensional very ample family of plane curves. By (3) and (4) in the definition, it will suffice to show that $\dim(L)=2$ and each $\dim(I_l)=1$. These follows from the ensuing three claims. Throughout, we will assume $(P,L,I)$ (and $X$) are $\emptyset$-definable.

\begin{claim}\label{C: pseudoplane lines have rank 1} For all $l\in L$ we have $\dim(I_l)=1$.
\end{claim}
\begin{claimproof} By assumption $I_l\subset X^2$ and $I_l$ is infinite, so we need only rule out the case that $\dim(I_l)=2$. But in this case $I_l$ is generic in $X^2$, so $\dim(X^2-I_l)=1$. By the finiteness of Morley degree, it follows easily (using (3)) that there can be only finitely many distinct $I_{l'}$ other than $I_l$. In other words we get that $L$ is finite, which clearly contradicts (2).
\end{claimproof}

\begin{claim} $\dim(L)\leq 2$.
\end{claim}
\begin{claimproof} Let $l\in L$ be generic, and let $p,p'$ be independent generics in $I_l$ over $l$. By (4), $\dim(l/pp')=0$, so $\dim(ppl')=\dim(pp')\leq 4$. But by the previous claim $\dim(pp'/l)=2$, so it follows that $\dim(l)\leq 2$, and thus $\dim(L)\leq 2$.
\end{claimproof}

\begin{claim} $\dim(L)\geq 2$.
\end{claim}
\begin{claimproof} Let $p\in P$ be generic, and let $l\in I^p$ be generic over $p$. By (2) we have $\dim(l/p)\geq 1$. Now if $\dim(L)<2$, we are forced to conclude that $\dim(l/p)=\dim(l)=1$. Thus $l$ and $p$ are independent. By assumption $\dim(p)=2$, so also $\dim(p/l)=2$. But then $\dim(I_l)\geq 2$, contradicting Claim \ref{C: pseudoplane lines have rank 1}.
\end{claimproof}
This completes the proof of the proposition. 
\end{proof}

\begin{remark} The above proof would be fairly straightforward, and included much sooner in the paper, if we knew right away that very ampleness gave us a two-dimensional very ample family of plane curves. The reason we had to wait until now to present the result is that we needed Proposition \ref{P: down} to get such a family.
\end{remark}

Before moving on, we point out one more fact that is relevant to pseudoplanes. As we have seen, very ampleness implies non-local modularity of strongly minimal structures. In Corollary \ref{C: 1-based iff no va type} below, we show, conversely, that non-local modularity (in fact, non-1-basedness) implies the existence of a very ample stationary type. This is true for stable theories, and the proof goes through \textit{complete type-definable pseudoplanes}. Recall the following (see \cite[\S, 4 Definition 1.6]{PillayBook}):

\begin{definition}\label{D: type pseudoplane} A \textit{complete-type-definable pseudoplane} consists of a complete type of a pair of (potentially imaginary) tuples $p=\operatorname{tp}(b,c)$ such that:
\begin{enumerate}
\item $b\notin\acl(c)$.
\item $c\notin\acl(b)$.
\item If $c'\neq c$ and $\operatorname{tp}(bc)=\operatorname{tp}(bc')$ then $b\in\acl(cc')$.
\item If $b'\neq b$ and $\operatorname{tp}(bc)=\operatorname{tp}(b'c)$ then $c\in\acl(bb')$.
\end{enumerate}
\end{definition}

The following is proved in (\cite[\S 4, Lemma 1.7]{PillayBook}):

\begin{fact}\label{F: 1-based iff no pseudoplane} A stable theory is 1-based if and only if it does not admit a complete-type-definable pseudoplane.
\end{fact}

\begin{remark} There is a subtlety in Definition \ref{D: type pseudoplane}. Namely, suppose $(P,L,I)$ is an arbitrary type-definable pseudoplane. Then $(P,L,I)$ might not contain any complete-type-definable pseudoplane: indeed, it is important that the partial type defining the incidence relation be complete. In fact, Fact \ref{F: 1-based iff no pseudoplane} fails if one considers general type-definable pseudoplanes. Let us sketch why. We leave it to the reader to verify that in a 1-based stable theory one can have a definable set $X$ equipped with a definable function $f:X\rightarrow X$ with all fibres infinite. One can then build a definable pseudoplane $(P,L,I)$ by letting $P=L=X$ and letting $I$ be two copies of $f$ (one in each direction).
\end{remark}

We now show:

\begin{corollary}\label{C: 1-based iff no va type} A stable theory admits a very ample non-algebraic stationary type if and only if it is not 1-based.
\end{corollary}
\begin{proof} Let $T$ be a stable theory. First suppose $T$ is 1-based. We show that no non-algebraic stationary type is very ample. To see this, let $p=\operatorname{tp}(a/c)$ be stationary, where $c=\operatorname{Cb}(p)$. Since $p$ is not algebraic, there is $b\models p$ with $b\neq a$. Now by 1-basedness, $c\in\acl(a)$, so automatically $c\ind_a b$, which shows that $p$ is not very ample.

Now suppose $T$ is not 1-based. By Fact \ref{F: 1-based iff no pseudoplane}, there is a complete-type-definable pseudoplane $\operatorname{tp}(bc)$. Let $d=\operatorname{Cb}(\stp(b/c))$. Clearly, $d\in\acl(c)$. On the other hand, by Definition \ref{D: type pseudoplane}, any $c$-conjugate of $\tp(b/c)$ is not parallel to $\tp(b/c)$. Thus $c\in\dcl(d)$, so in particular $c$ is interalgebraic with $d$.

We now claim that $\tp(b/d)$ is very ample. To see this, let $b'$ be any realization which is distinct from $b'$. We want to show that $b'$ forks with $d$ over $b$; by interalgebraicity, it is enough to show that $b'$ forks with $c$ over $b$. By symmetry, it is equivalent to show that $r=\tp(c/bb')$ is a forking extension of $q=\tp(c/b)$. But this is clear, since $r$ is algebraic and $q$ is not. 
\end{proof}

\section{Very Amplness in Strongly Minimal Groups}\label{S: applications}

In this final section we give some results related to very ampleness in strongly minimal groups. We first prove that non-locally modular strongly minimal groups admit non-affine plane curves (a result that has been long assumed but does not exist in writing). We then use this result to show that Question \ref{very ample sort question} has a positive answer for groups, and that divisible strongly minimal groups are already very ample; finally, we apply the very ampleness of divisible groups to to characterize fullness of strongly minimal ACF-relics with definable divisible group structures, in particular officially answering an old question on expansions of the multiplicative group.

\subsection{Existence of Non-Affine Plane Curves}

In \cite{HrPi87} Hrushovski and Pillay show that a stable group $\CG$ is 1-based if and only if every definable subset of $G^n$ (any $n$) is \emph{affine}, i.e., a boolean combination of cosets of definable subgroups of $G^n$. Restricted to strongly minimal (expansions of) groups, it has been well known among experts that non-1-basedness (equivalently, non-local modularity) even implies the existence of a non-affine plane curve. While this result is often used in practice, there does not seem to exist a full treatment anywhere in writing\footnote{The result is claimed  in \cite[Proposition 4.2]{KowRand}, but there is a mistake in the proof, that we do not see how to bridge.}. It can, however, be deduced by a particular application of Proposition \ref{P: large families}. Thus, we dedicate the present section to proving this result as a service to the community. 

Throughout this subsection, we fix $\CG$, a strongly minimal expansion of a group. We will say that a definable set $X\sub G^n$ is \emph{affine} if it is a (finite) boolean combination of cosets of definable subgroups of $G^n$; and \textit{almost affine} if it is almost equal to an affine set. Note that affine sets are almost affine, and the converse holds for one-dimensional sets (since finite sets are affine); more generally, note that arbitrary finite Boolean combinations of affine sets are affine. 

\begin{notation} For any non-empty definable set $X\subset G^n$, we let $\stab(X)$ denote the set of $g\in G^n$ such that $g+X$ almost coincides with $X$; so $\stab(X)$ is a subgroup of $G^n$ which is definable over any set of parameters defining $X$.
\end{notation}

We will show:

\begin{theorem}\label{T: nlm implies non affine} Let $\mathcal G=(G,+,\dots)$ be a strongly minimal expansion of a group. Then the following are equivalent:
\begin{enumerate}
    \item $\mathcal G$ is locally modular.
    \item Every definable plane curve in $\mathcal G$ is affine.
    \item Every definable one-dimensional set $X\subset G^n$, for any $n$, is affine.
    \item Every stationary definable set $X\subset G^n$, for any $n$, is almost affine.
    \item Every definable subset of each $G^n$ is affine.
\end{enumerate}
\end{theorem}

We will be brief, as many parts of the above proof are standard, and follow from either well known arguments or the main result of \cite{HrPi87}. Our  contribution will be the quick proof of the implication (3) $\Rightarrow$ (4) using Proposition \ref{P: large families}. In fact there is also a rather easy direct proof of the result of \cite{HrPi87} restricted to the strongly minimal case (using only the equivalence of local modularity and 1-basedness). So, for completeness, we now give this proof (so that our proof of Theorem \ref{T: nlm implies non affine} does not depend on \cite{HrPi87}). 

\begin{lemma}\label{L: sm HrPi}
    The following are equivalent:
    \begin{enumerate}
    \item $\CG$ is locally modular.
    \item If $\{X_t:t\in T\}$ is a non-empty faithful family of definable subsets of a definable set $Y$, each of dimension $k$, then $\dim(T)\leq\dim(Y)-k$.
    \item Every stationary definable set $X\subset G^n$ (for all $n$) is almost affine.
    \item Every definable set $X\subset G^n$ (for all $n$) is affine.
    \end{enumerate}
    Moreover, if (1)-(4) hold, then every definable subgroup $H\subset G^n$ (for all $n$) is $\acl(\0)$-definable.
\end{lemma}
\begin{proof}
    (1) $\Rightarrow$ (2): Assume $\CG$ is locally modular, and let $\{X_t:t\in T\}$, $Y$, and $k$ be as in (2). Without loss of generality $\{X_t\}$ and $Y$ are $\0$-definable. Let $t\in T$ be generic and $x\in X_t$ generic over $t$. Then $\dim(xt)=\dim(T)+k$, and by 1-basedness (equivalently local modularity) $t\in\acl(x)$, so $\dim(x)=\dim(T)+k$, thus $\dim(Y)\geq\dim(T)+k$, which is equivalent to (2).

    (2) $\Rightarrow$ (3): Let $X\subset G^n$ be stationary of dimension $k$, without loss of generality $\emptyset$-definable. Applying (2) shows that the family of translates of $X$ is at most $n-k$-dimensional (i.e. a generic translate of $X$ has canonical base of dimension at most $ n-k$). It follows that $\dim(\stab(X))\geq k$, which easily implies that $X$ is almost equal to a coset of $\stab(X)$, proving (3).

    (3) $\Rightarrow$ (4): By induction on dimension. Namely, assume (3) holds, let $X\subset G^n$ be definable, and assume that all definable sets of smaller dimension than $X$ are affine. By (3) applied to each stationary component of $X$, it follows that $X$ is almost affine. So there is an affine set $Y$ almost coinciding with $X$. Then $X$ is a Boolean combination of $Y$, $X\setminus Y$, and $Y\setminus X$, and the second two are of these are affine by induction.

    (4) $\Rightarrow$ (1): 
    Let $S$ be any strongly minimal plane curve; we show that $\dim(\operatorname{Cb}(S))\leq 1$, implying local modularity. Now we can write $S=C_{t_0}$ for some generic $t_0\in T$, where $\mathcal C=\{C_t:t\in T\}$ is a $\0$-definable (not necessarily faithful) family of plane curves, and $T\subset G^n$ for some $n$. Absorbing more constants to the language, if needed,  we may assume $T$ is stationary. It follows that the total space $C=\{(x,t):x\in C_t\}$ of $\mathcal C$ is stationary, so (by (4)) almost coincides with a coset $H$ of some connected definable subgroup of $G^{n+2}$. Note that (since $H$ is a coset) there is some (necessarily $\0$-definable) subgroup $K\leq G^2$ such that each non-empty fibre $H_t\subset G^2$ is a coset of $K$. In particular, since $t_0\in T$ is generic, it follows that $S=C_{t_0}$ almost coincides with a coset of $K$. Thus $\operatorname{Cb}(S)$ is definable over a single element of the one-dimensional group $G^2/K$, which gives $\dim(\operatorname{Cb}(S))\leq 1$, as desired.

    So we have shown that (1)-(4) are equivalent. Now suppose $\CG$ is locally modular. If there is a connected definable subgroup $H\subset G^n$ for some $n$ (say of dimension $k$) which is not $\acl(\0)$-definable, then there is a $\0$-definable infinite faithful family $\{H_t:t\in T\}$ of $k$-dimensional subgroups of $G^n$, whose generic members are connected. Passing to a subfamily and adding parameters if necessary, we may assume that $\dim(T)=1$. Then the family $\{X_s:s\in S\}$ of all cosets of all $H_t$ is a $\emptyset$-definable faithful $(n-k+1)$-dimensional family of $k$-dimensional subsets of $G^n$, which contradicts (2).
    \end{proof}
    
Before proceeding, let us also isolate the following two well known and easy facts: 
    
\begin{lemma}\label{L: plane curves to curves}
Let $\CH$ be a group of Morley rank 1 definable in an $\omega$-stable structure. Then any definable subset of $G$ is affine.  
\end{lemma}
\begin{proof}
Let $S\sub G$ be any definable set, and let $H_1, \dots, H_k$ be the cosets of $H^0$ (the definable connected component of $H$). Then each $H_i$ is strongly minimal, so $S\cap H_i$ is either finite or cofinite for all $i$. Thus $S$ has finite symmetric difference with $\bigcup \{H_i: \mr(S\cap H_i)=1\}$. Since the right-hand side is affine, set and the symmetric difference is finite (therefore also affine), the lemma follows. 
\end{proof}

\begin{lemma}\label{L: coset chunk} Let $S\subset G^n$ be stationary and $A$-definable. Then the following are equivalent:
\begin{enumerate} 
\item $S$ is almost affine.
\item For generic independent $x,y,z\in S$ over $A$, $x+y-z$ is also generic in $S$ over $A$.
\end{enumerate}
\end{lemma}
\begin{proof}
First assume (1) holds. Then by stationarity, $S$ almost coincides with a coset $H$ of a subgroup of $G^n$. Then $H$ clearly satisfies (2), so also $S$ does.

Now assume (2) holds. Fix $z\in S$ generic over the parameters defining $S$. Then (2), equivalently stated, gives that $S-z$ is almost contained in $\stab(S)$ (since given $x,y,z$ as in (2) it is clear that $(x-z)\in\stab(S)$). Thus $S$ is almost contained in $\stab(S)+z$. On the other hand, since $z$ is generic in $S$ it follows by definition that $\stab(S)+z$ is almost contained in $S$, so that in fact $S\sim\stab(S)+z$. Thus, (1) holds.
\end{proof}
    
Let us now proceed with the proof of Theorem \ref{T: nlm implies non affine}:

\begin{proof}
The implications (1) $\Rightarrow$ (2), (4) $\Rightarrow$ (5), and (5) $\Rightarrow$ (1) are contained in Lemma \ref{L: sm HrPi}. We show (2) $\Rightarrow$ (3) and (3) $\Rightarrow$ (4):

(2) $\Rightarrow$ (3): Let $C\subset G^n$ be definable of dimension 1. Since a finite union of affine sets is affine, we may assume $C$ is strongly minimal. Since $C$ is infinite, there is a projection $\pi:G^n\rightarrow G$ with cofinite image in $G$. Without loss of generality $\pi=\pi_1$ (the leftmost projection). Now for $i=2,\dots,n$, let $C_i=\pi_{1i}(C)\subset M^2$ (the image of $C$ in the first and $i$th coordinates). Then $C_i$ is a strongly minimal plane curve, so by (2) almost coincides with a strongly minimal coset $H_i$. Deleting finitely many points from $C$ if necessary, we may assume $C_i\subset H_i$ for all $i$.
    
Using that $H_i$ is a coset and $\pi_1(C)$ is cofinite in $G$, it follows that $\pi_1(H_i)=G$. So for each $i$, there is an element $c_i\in G$ with $(0,c_i)\in H_i$. Then replacing $C$ with $C-(0,c_2,\dots,c_n)$ if necessary, we may assume each $H_i$ is in fact a subgroup of $G^2$.
    
Now let $H=\{(x_1,\dots,x_n):(x_1,x_i)\in H_i\textrm{ for }i=2,\dots,n\}$. Then $C\subset H$, and it is easy to verify that $H$ is a one-dimensional definable subgroup of $G^n$. So by Lemma \ref{L: plane curves to curves}, $C$ is affine.

(3)$\Rightarrow$(4): If $\CG$ is locally modular this is immediate by Lemma \ref{L: sm HrPi}. So assume $\CG$ is not locally modular. Let $X\subset G^n$ be stationary of dimension $m$ and definable over $A$. Since finite sets are affine, we may assume that $m\geq 1$. Then by Proposition \ref{P: large families}, there is $B\supset A$ and a strongly minimal definable set $Y\subset X$ over $B$, such that any three independent generics in $Y$ over $B$ are also independent generics in $X$ over $A$. Let $x,y,z$ be as such. By (3) $Y$ is affine, so by Lemma \ref{L: coset chunk} $x+y-z$ is generic in $Y$ over $B$, and thus (again by the choice of $Y$) also generic in $X$ over $A$. By Lemma \ref{L: coset chunk} again, this implies that $X$ is almost affine.    
\end{proof}

\subsection{Very Ampleness in Groups}
We now apply Theorem \ref{T: nlm implies non affine} to show that strongly minimal groups have very ample finite quotients, and conclude that divisible strongly minimal groups are already very ample. We will need the following, which is an easy exercise (see \cite[Lemma 3.10]{ElHaPe} and Example \ref{group ex}): 

\begin{lemma}\label{L: stab char of va} Let $\mathcal G=(G,+,\dots)$ be a non-locally modular strongly minimal expansion of a group, and let $C\subset G^2$ be a strongly minimal plane curve. 
\begin{enumerate}
    \item The following are equivalent:
    \begin{enumerate}
    \item The family of translates $\{C+t:t\in T^2\}$ is faithful.
    \item The family of translates $\{C+t:t\in T^2\}$ is very ample.
    \item The group $\stab{(C)}$ is trivial.
\end{enumerate}
    \item The following are also equivalent:
    \begin{enumerate}
    \item $C$ is affine.
    \item $\stab(C)$ is infinite.
    \end{enumerate}
\end{enumerate}
\end{lemma}

Now we show:

\begin{theorem}\label{very ample group thm} Let $\mathcal G=(G,+,\dots)$ be a non-locally modular strongly minimal expansion of a group. Then:
\begin{enumerate}
    \item There is a finite subgroup $H\leq G$ such that $G/H$ is very ample. In particular, $\mathcal G$ admits a very ample sort.
    \item If $\mathcal G$ is divisible, then $\mathcal G$ itself is very ample.
\end{enumerate}
\end{theorem}
\begin{proof}
First we prove (1). By Theorem \ref{T: nlm implies non affine}, there is a strongly minimal plane curve $C\subset G^2$ which is not affine. Then by Lemma \ref{L: stab char of va}, the group $\stab{(C)}\leq G^2$ is finite. Let $H\leq G$ be a finite subgroup such that $\stab(C)\leq H^2$ (for example, $H$ could be the group generated by all coordinates of elements of $\stab(C)$). Let $C/H$ be the image of $C$ in $(G/H)^2$, by applying the projection $G\rightarrow G/H$ to both coordinates. It follows easily that, in the strongly minimal group $(G/H)$, the set $C/H$ is a strongly minimal plane curve with trivial stabilizer. So by Lemma \ref{L: stab char of va} again, $G/H$ is very ample.

Now assume $G$ is further divisible; we prove (2). Let $H$ be as in (1), and let $h=|H|$. Then the map $x\mapsto h\cdot x$ is well-defined on $G/H$, and so gives a definable map $G/H\rightarrow G$. Moreover, by divisibility this map is surjective with finite fibres. Thus, we are in the situation of Lemma \ref{L: very ample images} with $X=G/H$ and $Y=G$. So since $G/H$ is very ample, the lemma implies that so is $G$.
\end{proof}

\subsection{Expansions of Algebraic Groups}

 In the final subsection, we apply Theorem \ref{very ample group thm} to characterize ACF-definable expansions of 1-dimensional divisible algebraic groups. Our motivation is as follows. Suppose $(K,+,\times)$ is an algebraically closed field, and let $K^{\textrm{lin}}$ denote the structure $(K,+,\{\lambda_a:a\in K\})$, where $\lambda_a$ is the map $x\mapsto ax$. In the paper \cite{Martin}, Martin makes the following two conjectures:

\begin{enumerate}
    \item There are no intermediate structures between $(K,\times)$ and $(K,+,\times)$.
    \item If $\operatorname{char}(K)=0$ then there are no intermediate structures between $K^{\textrm{lin}}$ and $(K,+,\times)$. 
\end{enumerate}

(2) was proved by Marker and Pillay \cite{MaPi}; however, it seems that (1) has not since been addressed. Our result in this subsection will in fact be a more general statement implying both (1) and (2). Namely, we show:

\begin{theorem}\label{group trichotomy thm} Let $K$ be an algebraically closed field, and let $(G,\cdot)$ be the group of $K$-points of a one-dimensional divisible algebraic group over $K$. Let $G^{\mathrm{lin}}$ denote the structure endowing $G$ with the group operation and all of its endomorphisms (as an algebraic group), and let $G^{\mathrm{Zar}}$ denote the full $K$-induced structure on $G$. Then there are no intermediate structures between $G^{\mathrm{lin}}$ and $G^{\mathrm{Zar}}$.
\end{theorem}

\begin{remark} Since every endomorphism of the additive group is a scaling, and every endomorphism of the multiplicative group is a power map (thus definable from the group operation alone), we obtain (1) and (2) above. Moreover, note that Theorem \ref{group trichotomy thm} also applies to elliptic curves over $K$.
\end{remark}

\begin{proof} Let $\mathcal G=(G,\cdot,\dots)$ be a reduct of $G^{\mathrm{Zar}}$ which properly expands $G^{\mathrm{lin}}$.
\begin{claim} There is a non-affine $\mathcal G$-definable set $X\subset G^n$ for some $n$.
\end{claim}
\begin{proof} It suffices to show that every connected definable subgroup of $G^n$ is definable in $G^{\mathrm{lin}}$. So, let $H\leq G^n$ be a connected definable subgroup, say of dimension $d$. So $H$ is stationary, which implies there is an almost finite-to-one projection $H\rightarrow G^d$. Since $H$ is a subgroup, this implies that $H\rightarrow G^d$ is everywhere finite-to-one, and the fibers of $H\rightarrow G^d$ are cosets of a finite group, say $K\leq H$. Let $k=|K|$, and let $kH$ be the image of $H$ under scaling by $k$. Then the projection $kH\rightarrow G^d$ is the graph of a definable function $f:G^d\rightarrow G^{n-d}$. Since $kH$ is a subgroup, $f$ is a homomorphism. In particular, each coordinate component of $f$ (i.e. $G^d\rightarrow G^{n-d}\rightarrow G$) is a definable homomorphism from $G^d$ to $G$, and thus $f$ is the direct sum of $d$ definable endomorphisms of $G$. In particular, this implies that $f$ is $G^{\mathrm{lin}}$-definable, and thus so is $kH$.

Finally, let $H'$ be the preimage of $kH$ under scaling by $K$; so $H\leq H'$. Note that by divisibility, scaling by $k$ is finite-to-one on $G^n$; thus $\dim(H)=\dim(kH)=\dim(H')$. Then, since $H$ is connected, it must be the connected component $(H')^0$. In particular, $H$ is the image of $H'$ under scaling by $l$ for some $l$. Now to complete the proof, recall that $kH$ is $G^{\mathrm{lin}}$-definable; thus so is $H'$ by definition, and thus (scaling by $l$) so is $H$.
\end{proof}
Now by the claim and the main result of \cite{HrPi87} (or Lemma \ref{L: sm HrPi} if desired), it follows that $\mathcal G$ is not locally modular. Then by divisibility and Theorem \ref{very ample group thm}, $\mathcal G$ is very ample. Finally, by \cite{HaSu} $\mathcal G$ satisfies the Zilber trichotomy; so by Theorem \ref{very ample acf char thm}, $\mathcal G$ is full in $K$. In other words, $\mathcal G$ is interdefinable with $G^{\mathrm{Zar}}$, as desired.
\end{proof}

\end{document}